\documentclass[12]{amsart}
\usepackage{amsmath}
\usepackage{amssymb}
\usepackage{mathrsfs}
\usepackage{stmaryrd}

\usepackage[all]{xypic}
\usepackage{amsfonts}
\newtheorem{theorem}{Theorem}[subsection]
\newtheorem{lemma}[theorem]{Lemma}
\newtheorem{prop}[theorem]{Proposition}
\newtheorem{co}[theorem]{Corollary}
\theoremstyle{definition}
\newtheorem{definition}[theorem]{Definition}
\newtheorem{example}[theorem]{Example}
\newtheorem{remark}[theorem]{Remark}

\theoremstyle{remark}

\newtheorem{con}[theorem]{Convention}
\numberwithin{equation}{subsection}

\def \m {\mathfrak m}
\def \G {\mathcal G}
\def \H{\mathcal H}
\def \bt {Barsotti-Tate\ }
 
\def \K {\mathcal K}
\def \Z {\mathbb Z}
\def \inj {\hookrightarrow }
\def \to {\rightarrow}
\def \spec \text{spec}
 
\def \tensor {\otimes }
\def \simto {\simeq}\title{Research Statement }

\def \e {e}

\def \ab {\mathscr A}

\DeclareMathOperator{\gal}{Gal}

\def \Q {\mathbb Q}
\def \t {\textnormal}
\def \Z {\mathbb Z}
\def \w {\widetilde}
\def \K {\mathcal K }
\def \O {\mathcal O}
\def \ito {\overset  \sim  \to}
\def \D {\mathbb D}
\begin{document}

\title{Potentially Good Reduction of Barsotti-Tate  Groups}

\author{Tong Liu}
\email{Tong.Liu@math.u-psud.fr}
\address{Universit\'e de Paris-Sud,Math\'ematiques, B\^{a}timent
425, 91405 Orsay Cedex  France}







\begin{abstract}
 Let $R$ be a complete discrete valuation ring of mixed characteristic (0, $p$) with
 perfect residue field, $K$  the fraction field of $R$.
 Suppose $G$ is a \bt  group ($p$-divisible group) defined over $K$ which acquires good reduction over a finite extension $K'$ of $K$. We prove that
 there exists a constant $c\geq 2 $ which depends on the absolute ramification index $e(K'/\mathbb Q_p)$ and the height of $G$ such that
  $G$ has good reduction over $K$ if and only if $G[p^c]$ can be extended to a finite flat group scheme over $R$. For abelian varieties with
  potentially good reduction, this result generalizes  Grothendieck's ``$p$-adic N\'eron-Ogg-Shafarevich criterion" to finite level.
  We use methods that can be generalized to study semi-stable $p$-adic  Galois representations with general
Hodge-Tate weights, and in particular leads to a proof of a
conjecture of Fontaine and gives a constant $c$ as above that is
independent of the height of $G$.
\end{abstract}

\maketitle


\tableofcontents
\section{Introduction}
Let $R$ be a complete discrete valuation ring of mixed characteristic (0, $p$), $K$ the fraction field, $\m $ the maximal ideal of $R$,   $k=R/\m$
 the residue field and $A$ an abelian variety defined over $K$. We say $A$ has \emph{good reduction}
if there exists an abelian scheme $\ab$ over $R$ such that $ \ab \otimes_R K\simeq A$.
 We say  $A$ has \emph{potentially good reduction} if there exists a finite extension $K'/K $ and an abelian scheme
$\ab'$ over $R'$ such that $\ab' \otimes_{R'} {K'} \simeq A\otimes_K {K'}$.
 Let $I $ be the inertia subgroup of
${\rm {Gal}}(\overline K/K)$ and $l\neq p $ be any prime. Using the
criterion of N\'eron-Ogg-Shafarevich \cite{serre}, it is not
difficult to prove that if $A $ has potentially good reduction then
$A$ has good reduction if and only if the action of $I $ on $A[l^c]$
is trivial, where $c=1$ for $l\neq 2$ and $c=2$ for $l=2$.

 When $l=p$, the problem is more subtle. For example,  the relevant finite flat group schemes are not  \'etale, and thus
may not be determined by
their generic fiber as in the case $l\neq p $.  Work of Grothendieck \cite{sga7} and Raynaud~\cite{ray} proves that $A$ has good reduction if and only if
for all $n$,  $A[p^n]$ can be extended to a finite flat group scheme over $R$. In particular, the property of good reduction (and hence
potentially good reduction) is encoded in the \bt group $G=\varinjlim A[p^n]$
(also known as the $p$-divisible group associated to $A$). This fact
is implicit in the work of Wiles and others on modularity of certain
Galois representations of weight 2 (\cite{bcdt},~\cite{tw}). In what
follows, we fix $l=p$ and study potentially good reduction for \bt
groups.

Let $K'$ be a finite  extension of $K$, $R' $ the ring of integers
of $K'$. Suppose a \bt group $G$ with height $h$ defined over $K$
has {\it good reduction} over $K'$ in the sense that there exists a
\bt group $\G'$ over $R'$ such that $\G'\otimes_{R'} K'\simeq
G\otimes_K K' $. By \S2, \cite{ray}, $G$ has good reduction over $K$
if only if $G[p^n]$ can be extended to a finite flat group scheme
over $R$ for all $n$. The following question was raised  by N. Katz:

{\textit{For $G$ as above with height $h$ and good reduction over
$K'$, does there exist a constant $c=c(K',K, h)$ such that $G$ has
good reduction over $K$ if and only if $G[p^c]$ can be extended to a
finite flat group scheme over $R$? If so, how does $c$ depend on the
arithmetic of $K'/K$, and how does it depend on  the height $h$? } }

The general case is easily reduced to the case when $K'/K$ is totally ramified. Let $\e $ be the absolute ramification
index of $K'/\mathbb Q_p$.  In the case $\e< p-1$,
 B. Conrad  solved this question in \cite{bd}.
He proved that if $G[p]$ can be extended to a finite flat group scheme, then $G$ has
good reduction. That is,  when $\e< p-1$, $c$ can be taken to be 1 for all $h$. His idea was to compute descent data for the generalized
 Honda system associated to $G$. Unfortunately, in the case $\e\geq p-1$, there is no Honda system available. Although the theory of Breuil modules \cite{b1}
 is available  in this case (at least for $p>2$), its behavior under base change is very complicated due to
 the dependence of theory on the choice of uniformizer, so it is hard to apply Breuil's theory for this question (especially  when $p\mid \e$).

We use Messing's deformation theory \cite{mess} to attack this
question and get the following result:
\begin{theorem}[Main Theorem]\label{main}
Assume that $k$ is perfect. Suppose $G$ is a \bt  group over $K$
that acquires good reduction over a finite extension $K'$. There
exists a constant $c\geq 2 $ which depends on $\e=e(K'/\mathbb Q_p)$
and the height of $G$ such that
  $G$ has good reduction over $R$ if and only if $G[p^c]$
can be extended to a finite flat group scheme over $R$.
\end{theorem}

As we will explain in Remark \ref{remark2}, a refinement of the
methods of this paper allows us to take $c$ to be independent of the
height of $G$ (though in practice one usually studies $G$'s with a
fixed, or at least bounded, height). Applying this result to the
case of abelian varieties and using the precise form of the
semi-stable reduction theorem for abelian varieties (i.e., bounding
the degree of the extension over which semi-stable reduction is
attained), we get :

\begin{theorem}\label{th1.2}
Assume that $k$ is perfect. Suppose $A $ is an abelian variety
over $K$ with potentially good reduction. There exists a constant
$c'\geq 2 $ depending
 on $e(K/\mathbb Q_p)$ and $\dim A$
 such that   $A$ has good reduction over $K$ if and only if $A[p^{c'}]$ can be extended to a finite flat group scheme over $R$.
\end{theorem}

 We may and do assume $R'/R$ is totally ramified and generically Galois.
 Let $\Gamma= \gal(K'/K)$ and $G' =G\otimes_K K'$.  Since $G'$ has good reduction over $K'$, there
 exists a \bt group $\G'$ over $R'$ such that $\G'\otimes_{R'}K'\simeq G'$. Using the $K$-descent $G$ of $G'$, clearly
 $G'$ has natural descent data $\phi _{\sigma }$, where $\phi_ {\sigma } : \sigma^* (G') \simeq G'$ is a collection of
 $K'$-automorphisms satisfying
the cocycle condition. By Tate's Theorem \cite{ta}, the descent data
on the generic fiber will uniquely extend over $R'$. That is, we get
isomorphisms $ \widetilde \phi_\sigma: \sigma^*(\G')\simeq \G'$ over
$R'$ satisfying the cocycle condition, but this is \emph{not} fppf
descent data with respect to $R\to R'$ when $R\to R'$ is not
\'etale.
 It is obvious  that $G$ has good reduction if and only if
 there exists a \bt group $\mathcal  G  $ over $R$ and an $R'$-isomorphism
$\mathcal  G \otimes  _R  R' \simeq  \mathcal  G'$ compatible with
the Galois action (using $\widetilde\phi_\sigma$'s on $\G'$). In
general, given a finite Galois group $\Gamma$,  we call a pair
$(\G',\widetilde\phi_ {\sigma }  )$ a \emph{$\Gamma$-descent datum}
where  $\G'$ is a finite flat group scheme (a \bt group) over $R'$
and  $ \widetilde \phi_\sigma: \sigma^*(\G')\simeq \G'$  are
$R'$-isomorphisms satisfying the cocycle condition. We say a
$\Gamma$-descent datum  $(\mathcal  G ' , \widetilde\phi_ {\sigma }
)$ has an \emph{effective descent} if there exists a finite flat
group scheme (a \bt group) $\G $ over $R$ and an $R'$-isomorphism
$\mathcal  G \otimes  _R  R' \simeq  \mathcal  G'$ compatible with
the Galois action (using $\widetilde\phi_\sigma$'s on $\G'$,  see \S
2 for more details).
 Using deformation theory of \bt groups,
\S2 to \S 4 will prove the following:

 \begin{theorem}[Descent  Theorem]\label{th1}
The data
$(\mathcal  G ' , \widetilde\phi_ {\sigma } )$ has an effective descent to a \bt group over $R$ if and only if
 $(\mathcal  G '[p^{c_0}], \widetilde\phi_ {\sigma })$ has an effective descent to a finite flat group scheme over $R$, where $c_0\geq 1$ only depends on $\e$.

\end{theorem}
\begin{remark}We will prove the same  Theorem  in the equi-characteristic case. In such a case, the constant
$c_0 $  depends on  the \emph{relative} discriminant $
\Delta_{K'/K}$ and the height of $G$. We do not optimize the
constant $c_0$, so there is still room to improve.
\end{remark}
The Descent Theorem is the heart part of this paper. Let us sketch
the idea of proof as the following: In \S 2, we develop  an \emph{ad
hoc} finite Galois descent theory on torsion level (e.g.,
$R/\m^n$-level). Using such language, $G$ has good reduction over
$K$ if and only if for each $n\geq 1$, the $\Gamma$-descent datum
$(\G'_n, \widetilde\phi_{\sigma, n })$ has an effective descent to
$\G_n$ over $R/\m ^n $ and compatible with base change
$R'/\m^{n+1}\to R'/\m^{n}$, where $\G'_n =\G' \otimes_{R'}
R'/\m^{n}$ and  $\widetilde\phi_{\sigma, n }=\widetilde\phi_\sigma
\otimes_{R'} R'/\m^{n}$. Since $R'$ is totally ramified over $R$, we
see that if such $\G_n $ does exist, $\G_n$  and $\G'_n$ must has
same closed fiber $\G_0$. Thus, $\G_n$ can be seen as some
``correct" deformations of $\G_0$. Luckily, Messing's theory
provides some crystals to precisely describe deformations of \bt
groups. So the proof of Theorem \ref{th1} can be reduced to the
effectiveness of descending  a short exact sequence of finite free
$R'$-modules with semi-linear $\Gamma$-actions coming from $G'$:
$$0\to M_1'\to M_2' \to M_3' \to 0.$$
In general, unlike the ordinary Galois descent theory, the effective
descent of such sequence does not necessarily exist. But work of \S
2 will prove that such effective descent does exist if for each
$i=1,2,3 $, $M_i'\otimes _{R'} R'/\m ^d  $ has an effective descent
to an $R/\m^d$-module $M_i$ for some constant $d\geq 1$ which
depends only on $e$ in the mixed characteristic case and depends
only on $\Delta_{K'/K}$ and the height of $G$ in the
equi-characteristic case.  Finally, in \S 3, we use deformation
theories from Grothendieck and Drinfeld to show  that  the effective
descents of $M'_i\otimes_{R'} R' /\m^d $ does exist by proving that
$\G'_d$ has an effective descent to $ \G_d$ under the hypothesis
that $\G'[p^{c_0}]$ has an effective descent. Then Theorem \ref{th1}
follows.

Note that Theorem \ref{th1} does not immediately prove Theorem 1.1.
In fact, although $G[p^{c_0}]$ can be extended to  a finite flat
group scheme $\G_{c_0}$ over $R$, $\G_{c_0}$ is not necessarily an
effective descent of  $\G'[p^{c_0}]$ because the extensions of a
finite $K$-group to a finite flat $R$-group scheme are not
necessarily unique. Hence, in \S5, we generalize Tate's isogeny
theorem \cite{ta} and Raynaud's result \S 3, \cite{ray} to finite
level:

\begin{theorem}\label{th2}
Let $\G $ be a truncated \bt group over $R$ of level $n$,
 $\H$ a finite flat group scheme over $R$ and $f_K: \G_K \to \H_K $ a $K$-group scheme morphism.
 Then  there exists a  constant $c_1\geq 2$ depending on the absolute ramification index $e(K/\Q_p)$ and
the height of $\G_K$ such that for all  $n\geq c_1 $
there exists an $R$-group scheme morphism $F : \G \to \H $ satisfying  $F  \otimes_R K = p^{c_1}\circ f_K $.
\end{theorem}

\begin{prop}\label{th3}
 Let $\G$  be a finite flat group scheme over $R$ whose generic fiber $\G_K$ is a truncated
\bt group of level $n$.
There exists a constant $c_2 \geq 2$ depending on $e$ and the height of $\G_K$
such that if  $n\geq c_2$
  then there exists a truncated \bt group $\G' $over $R$ of level $n-c_2$
and $R$-group scheme morphisms
$g : \G \to \G'$ and $g':  \G'\to \G  $ satisfying
\begin{enumerate}
\item $g'_K = g'\otimes _R K $ factors though $\G_K[p^{c_2}]$ and $g'_K : \G'_K \to \G_K[p^{c_2}]$ is an isomorphism
\item $(g\circ g' )\otimes _R K =p^ {c_2}$
\end{enumerate}
\end{prop}

Combining Theorem \ref{th2} and Proposition \ref{th3}, in the
beginning of \S 5, we will see that  $c=c_0+c_1 +c_2$,  and if
$G[p^c]$ can extend a finite flat group scheme $\G_c$ then there
exits a subgroup scheme $\G_{c_0}$ of $\G_c$ such that $\G_{c_0}$
provides an effective descent of $\G'[p^{c_0}]$, hence prove the
main Theorem.

\begin{remark}\label{remark2} Using  the
classification of finite flat group schemes over $R$ in \S 2.3,
\cite{kisin2}, we can prove that  $c_1$ and $c_2$, hence $c$  can be
chosen to depend only on the absolute ramification index and not on
the height. In fact, we can prove Theorem \ref{th2} and Proposition
\ref{th3} in the more general setting for \emph{$\varphi$-modules of
finite height} in the sense of Fontaine (\cite{fo4}). These results
will play a central technical role  in our proof of a conjecture of
Fontaine (\cite{fo10}) concerning that the semi-stability for
$p$-adic Galois representations is preserved under suitable inverse
limit processes, subject to a natural boundedness hypotheses on the
Hodge-Tate weights. The methods of this paper may be viewed as a
simpler version of the methods we develop in \cite{liu2} where we
prove Fontaine's conjecture.
\end{remark}

When writing of this paper was nearly  complete, Bondarko posted his
preprint \cite{Bond} on the {\tt arxiv}, where he gave a new
classification of finite flat group schemes over $R$ and obtained a
similar result (c.f. Theorem 4.3.2) as our main result, Theorem
\ref{main} above. We first remark that the hypotheses in our Main
Theorem are  weaker than that of Theorem 4.3.2 in Bondarko's
preprint \cite{Bond}. One of his hypotheses requires the finite flat
group scheme which extends the given generic fiber has to have the
correct ``tangent space". But our hypothesis does not have such a
restriction on such finite flat group schemes. Also here we use a
totally different method. Bondarko's theory mainly depends on the
theory of Cartier modules that gives a classification of formal
groups. His method works without a perfectness hypothesis on the
residue field and gives slightly better constants (cf. Theorem
4.1.2). But our approach by deformation theory can be adapted to
apply to rather general semi-stable $p$-adic Galois representations
without restricting to crystalline representations with Hodge-Tate
weights in $\{0,1\}$ (see \cite{b1}).  In contrast, Bondarko's
method of formal groups seems difficult to apply to such a wider
class of representations. To be more precise, the problem we solve
in this paper is a special case of the following more general
question:

\emph{Let $T$ be a $\gal(\overline K/ K)$-stable $\Z_p$-lattice in a
potentially semi-stable $p$-adic Galois representation $V$. Does
there exist a constant $c$ such that if $T/p^cT$ is torsion
semi-stable then $V$ is semi-stable? If so, how does $c$ depend on
the arithmetic of $K$ and $V$?}

Using the machinery provided by (\cite{kisin2}), the author has made
some progress on this question by modifying the method in this
paper. Also an extended version of Theorem \ref{th2} and Proposition
\ref{th3} is  used in our proof of  a conjecture of Fontaine
(\cite{liu2}), as explained in Remark \ref{remark2}.

\begin{con} From \S 2 to \S 4, we only assume $R$ is a complete discrete valuation ring with
perfect residue field $k$  with char$(k)=p>0$,  and we select a
uniformizer $\pi $ (resp. $\pi'$) in the maximal ideal $\m $
(resp. $\m'$). In particular, we do not assume that the fraction
field $K$ of $R$ has characteristic 0. We fixed a finite Galois
extension $K'/K $ with Galois group $\Gamma$. We have to consider
torsion modules such as  $M/\pi^n M$. In order to simplify  the
notation, we write $M/\pi^n $ to denote $M/ \pi^n M $. We use
``$\t{Id}$'' to denote the identity map (of modules, schemes,
etc.). If $M $ is an $R'$-module equipped with semi-linear action
of $\Gamma$, then $M^{\Gamma}$ denotes the $R$-submodule of
$\Gamma$-invariant elements.
\end{con}

{\bf Acknowledgment}. I thank my  advisor Brian Conrad for his
patient guidance and advice. Also I would like to thank A.J. de Jong
and James Parson for useful comments.

\section{Elementary Finite Descent Theory }
\subsection{Preliminary}
We will first set up an \emph{ad hoc}  ``descent theory"   from
$R'$ to $R$ using $\Gamma$-actions. We will see that for ``descent
data" on  a finite free $R'$-module, the obstruction to effective
descent can be determined by working  modulo $\m^n$  for a
sufficiently large $n$ depending only on the arithmetic of $K'/K$.

\begin {definition}  Let $M'$, $N'$ be $R'$-modules. Fix $\sigma \in \Gamma$.
A  $\sigma$-\emph{semi-linear homomorphism} $\phi_\sigma :
M'\rightarrow N' $ is  an  $R$-linear homomorphism such that for
all $x\in R' $ and    for all $m\in M' $,
$\phi_\sigma(xm)=\sigma(x)\phi(m) $. An $R'$-\emph{semi-linear
homomorphism} from $M'$ to $N'$ is a collection
 $(\phi_\sigma)_{\sigma \in \Gamma}$ where each $\phi_\sigma$ is $\sigma $-semi-linear and $\phi_{\sigma \tau}= \phi_\sigma \circ \phi_\tau $ for all $\sigma , \tau \in \Gamma$.
\end{definition}

It is obvious that the set of all $R'$-semi-linear automorphisms
of $M'$ forms a group, denoted by ${\text{Aut}}_{\Gamma}(M')$

\begin{definition} A  \emph{$\Gamma$-descent module} structure on  an  $R'$-module $M'$   is a pair
$(M', \phi )$,  where $\phi : \Gamma\to {\text{Aut}}_{\Gamma}(M')$
gives an $R'$-semi-linear action of $\Gamma$ on $M'$.
\end{definition}

Let $(M', \phi_\sigma )$, $ (N', \psi_\sigma)$ be two
$\Gamma$-descent modules. A \emph{morphism} $f: (M', \phi_\sigma )
\to (N', \psi_\sigma)$ is a $R'$-linear morphism $f: M'\to N'$
such that for all $\sigma \in \Gamma $, $f\circ
\phi_\sigma=\psi_\sigma \circ f $.  That is, a morphism is just an
$R'$-module morphism compatible with the $\Gamma$-actions.
  Most of time we use terminology ``Galois-equivariantly" instead of ``morphism of
$\Gamma$-descent modules".

\begin{example}\label{eg1}
Let $M$ be an $R$-module, so $M'= M\otimes_R R'$ has a natural
$\Gamma$-descent datum $\phi_\sigma$ that comes from the Galois
action of $\Gamma$   on  $R'$.
 We always omit $\phi_\sigma $ when we mention the canonical $\Gamma$-descent module $(M\otimes _R R', \phi_\sigma )$.
\end{example}
\begin{definition} A $\Gamma$-descent module $(M', \phi_\sigma)$ has an \emph{effective descent}   if there exists an $R$-module $M$ such that
$ M\otimes_R R'$ is isomorphic to $ ( M', \phi_\sigma)  $ as
$\Gamma$-descent modules,  or equivalently, if there  is  a
 $R'$-isomorphism $i: M\otimes_R R ' \simeq
M'$ such that $i$ is Galois equivariant.
\end{definition}

 Similarly, we can  define exact sequences of $\Gamma$-descent modules and effective  descent for such exact sequences.
 Also we can generalize such  descent language to the setting of schemes. For example, we can define
 $\Gamma$-descent data of an affine scheme (affine group scheme) over $R'$ by defining $\Gamma$-descent data of
its coordinate ring (Hopf algebra). Similar definitions  apply to
 the notion  of effective descent of an affine scheme (affine group scheme) over $R'$. More precisely, let
$X'= \t{Spec}(A')$ where $A'$ is an $R'$-algebra. It is easy to
see that a $\Gamma$-descent data structure on $A'$ is equivalent
to a collection of $R'$-isomorphisms $\Phi_\sigma:
\sigma^*(X')\simeq  X'$ satisfying cocycle conditions
$\Phi_{\sigma \tau }= \Phi_{\sigma} \sigma^*(\Phi_\tau)$ for all
$\sigma, \tau \in \Gamma$.   Thus we also call such a pair $(X',
\Phi_\sigma)$ a \emph{$\Gamma$-descent datum}. Let $(X'',
\Phi'_\sigma)$ be another $\Gamma$-descent datum.  A
\emph{morphism} between $\Gamma$-descent data $f : (X',
\Phi_\sigma)\to (X'', \Phi'_\sigma)$ is morphism of $R'$-schemes
such that the following diagram commutes for all $\sigma \in
\Gamma$:
\begin{equation} \label{diag2.1}
\begin{split}
\xymatrix{\sigma^*(X')\ar[d]_{\Phi_\sigma}^\wr \ar[r]^{\sigma^*(f)} & \sigma^*(X'')\ar[d]^{\Phi'_\sigma}_\wr\\
X' \ar[r]^f & X'' }
\end{split}
\end{equation}

Similarly, we define $f$ to be an  \emph{isomorphism} of
$\Gamma$-descent data if $f$ is an isomorphism of underlying
$R'$-schemes. The inverse is obviously a morphism of
$\Gamma$-descent data.
\begin{example}
Let $X= \t{Spec}(A)$ be an  $R$-scheme. Then we have a canonical
$\Gamma$-descent data $\iota _\sigma $ on $X\otimes_R R'$ coming
from the canonical $\Gamma$-descent algebra $A \otimes_R R'$.
\end{example}

We say that a $\Gamma$-descent datum $(X', \Phi_\sigma)$ has an
\emph{effective descent} if there exists an $R$-scheme $X$ such
that $(X\otimes_R R', \iota_ \sigma)$ is isomorphic to $(X',
\Phi_\sigma)$ as $\Gamma$-descent data. In the affine case, it  is
equivalent to say that the coordinate ring of $X'$ has an
effective descent as a $\Gamma$-descent algebra.
\begin{example}
As in \S1, let $G$ be a \bt group over $K$ which acquire good
reduction over $K'$ in the sense that there exists an \bt group
$\G'$. By Tate's Theorem \cite{ta} (De Jong's Theorem
\cite{dejong} in the equi-characteristic case),  we have
a $\Gamma$-descent datum $(\G', \w \phi_\sigma)$  coming from the
generic fiber $G \otimes_K K'$. Furthermore, $G$ has good
reduction over $K$ if and only if $(\G', \w \phi_\sigma)$ has an
effective descent and Theorem \ref{th1} claims that such effective
descent does exist if and only if $(\G[p^{c_0}], \w \phi_\sigma)$
has one for some constant $c_0$.
\end{example}
Let $(X' , \Phi_\sigma)$ and $(X'', \Phi'_\sigma)$ be
$\Gamma$-descent data and $f : X' \ito X'' $  an isomorphism  of
$R'$-schemes. The following easy lemma will be very useful to
check  if $f$ is an isomorphism of $\Gamma$-descent data.
\begin{lemma} \label{check} Using notations as above, for all $\sigma \in \Gamma$, define
$$
f_\sigma: \xymatrix{X' \ar[r]_-{\sim } ^-{\Phi^{-1}_\sigma}&
\sigma^*(X') \ar[r]^-{\sigma ^* (f)} & \sigma^*  (X'')
\ar[r]^-{\Phi'_\sigma}_-{\sim}& X''}.
  $$
Then $f$ is an isomorphism of $\Gamma$-descent data if and only if
$f^{-1}\circ f_\sigma = \t{Id}$ for all  $\sigma \in \Gamma$.
\end{lemma}
\begin{proof} From the diagram \eqref{diag2.1}, this is clear.
\end{proof}
\subsection{``Weak" full faithfullness} Define a functor  from the
category of $R$-modules to the category of $\Gamma$-descent
modules  over $R'$:
$$F: M\longmapsto (M\otimes_R R', \phi_\sigma)$$
as in Example \ref{eg1}. Unlike  ordinary descent theory (e.g.,
{\it fppf} descent theory), the functor $F$ is neither full nor
faithful.

\begin{example} Let $I '\subset R'$ be an ideal not arising from $R$. The $\Gamma $-action on $R'$ preserves $I'$, and
$I: = I '^\Gamma $ is a nonzero principal ideal of $R$, but the
map $I \otimes _R R' \inj I' $ is not surjective. That is, $I'$
with its evident $ \Gamma $-descent datum does not have an
effective descent.
\end{example}

\begin{example}If $(M', \phi_\sigma)$ is a $\Gamma$-descent module,  then
$\overline { M'}=M' /\pi^n M '$ has a canonical descent datum
$(\overline{M'},  \bar \phi_\sigma) $. We still use $\phi_\sigma$
to denote the semi-linear Galois action on $\overline {M'}$.
Similarly, if $(X, \Phi_\sigma)$ is a $\Gamma$-descent data, we
have natural  $\Gamma$-descent data $\Phi_\sigma \otimes_ R R/ \pi
^n $ on $X \otimes_R R/\pi^n$.
\end{example}

\begin{example} \label{example}
 Let  $M'= R'/ {\pi^n}$. There is a natural descent datum
$\phi_\sigma $ on $M'$  from the Galois action of $\Gamma $ on
$R'$. We see  that $R/\pi^{n}\to   (R'/ {\pi^n})^{\Gamma}$ is
injective. However in general, the equality does  not necessary
hold (In fact, $\t{H}^1(\Gamma, R')[\pi^n] $ is not necessarily
trivial, especially when $p|e(K'/K)$). Select an $\alpha \in
(R'/\pi^n)^\Gamma- R/\pi^n$ such that $R'/\pi^n =(R'/\pi^n)\cdot
\alpha$.
 Let $M_1=R/\pi^n $ and $M_2= (R/\pi^n)\cdot  \alpha $,
 so $M_1 $ and  $M_2$ are two  effective descents of $(M', \phi_\sigma) $. Note that the identity $\t{Id}: M'\to M'$ is a morphism between
 these two descent data on $M'$, but there does not exist a morphism  $f: M_1 \to M_2  $ such that $f\otimes_R R' =\t{Id}. $
\end{example}

As we have seen above, the  functor $F$ is not a fully faithful
functor.  However, we will show $F$ actually has ``weak" full
faithfulness in the following sense.

\begin{prop}{\label {func}}
Let $f: (M'_1, \phi_1 )\to (M'_2, \phi_2)$ be a morphism of
$\Gamma$-descent modules with  $M'_1 $, $M'_2$ finite free
$R'/\pi^n$-modules. Suppose $M_1 $, $M_2$ are effective descents
for $ M_1' $, $ M'_2 $ respectively. Then $M_1$ and $M_2$ are
finite free $R/\pi^n$-modules. Furthermore, there exists a
constant $r$ depending on $R'$ and $R$ such that for all $n\geq r$
the morphism
 $\pi^{r}f$  has a unique effective descent to a morphism $\hat f: M_1 \to M_2$.
\end{prop}

In  Lemma \ref{bound} below, we will see how the constant $r$
depends on the arithmetic of $R'$ and $R$.
\begin{co}
If $f$ in  Proposition \t{\ref{func}} is an isomorphism, then the
natural map $\bar f : M_1/\pi^{n-r}\to M_2/\pi^{n-r}$ induced by
$\hat f $ is an isomorphism.
\end{co}
\begin{co}\label{free}
 With notations as above,  if  $M_1'$ is a finite free $R'$-module then the effective descent of $(M_1,  \phi_1)$ is unique up to the unique isomorphism if it exists.
  Moreover, the formation of such descent \t{(}when it exists\t{)} is uniquely effective on the level of morphisms.
\end{co}
\begin{remark}
Of course, Corollary \ref{free} is obvious from considering $M_1'
\otimes_{R'} K'$ over $K'$. However for our purposes what matters
is that such an effective descent can be studied
 using   the effective descent at  torsion level.
\end{remark}

Let  $S':\  0 \to (M'_1, \phi_1) \to (M'_2, \phi_2)\to (M'_3,
\phi_3)\to 0$  be an exact sequence of $\Gamma$-descent modules.
Recall that we say that $S'$ has an effective descent  if there
exists an \emph{exact} sequence of $R$-modules  $S: \  0\to M_1
\to M_2 \to M_3 \to 0 $ such that the following diagram commutes.
$$
\xymatrix{0 \ar[r]& M_1' \ar[r] &M_2'\ar[r] & M_3' \ar[r] & 0\\
0 \ar[r] & M_1 \otimes_R R'\ar[r] \ar[u]^{f_1}_\wr& M_2 \otimes_R
R' \ar[r] \ar[u]_\wr^{f_2}& M_3 \otimes _R R'
\ar[r]\ar[u]^{f_3}_\wr & 0  }
$$
where $f_i$ is a Galois-equivariant  isomorphism for $i=1,2, 3$.
Using Corollary \ref{free}, we get:

\begin{co}\label{exact} Suppose $M_i'$ is a  finite free $R'$-module for each $i=1,2,3$. Then
$S'$ has an effective descent if and only if $(M'_i, \phi_i)$ has
an effective descent for $i=1,2,3$.
\end{co}

Before proving Proposition \ref{func}, we  need a  lemma to bound
Galois cohomology and precisely describe the constant $r$ in
Proposition \ref{func}. Let $\Delta _{K'/K}$ be the discriminant
of $R'/R$,  $r_0$  the $p$-index of $\Gamma$ (that is
$p^{r_0}\parallel  \#(\Gamma)$).
 We normalize the valuation $v(\cdot)$ on $R'$ by putting  $v(\pi)=1$, so  $v(K'^\times )=e(K'/K)^{-1}\cdot\Z$ where
 $e(K'/K)$ is the relative ramification index.
Finally, we put $r= v(\Delta_{K'/K})$ in the equi-characteristic
case and $r= r_0v(p)$ in the mixed characteristic case.

\begin{lemma}\label{bound} Let $(M', \phi_\sigma)$ be a $\Gamma$-descent module with
$M'$ a finite free $R'$-module of rank $m$.  The element
 $\pi^ r $ kills $\t{H}^1(\Gamma, M ')$ in  the mixed characteristic case. In the equi-characteristic case,
$\t{H}^1(\Gamma,  M')=0 $ if $r_0=0$ and $\pi^{m(r+1)}$ kills
$\t{H}^1 (\Gamma, M')$ if $r_0 >0$.
\end{lemma}
\begin{proof}
 By Hilbert's Theorem 90, $\t{H}^1(\Gamma, M'\otimes_{R'}K')=0$. Thus,   there exists a positive integer $\tilde r\geq 0$
such that $\pi ^ {\tilde r}$ kills the finite $R$-module
$\t{H}^1(\Gamma, M' )$. Hence, in both cases (equi-characteristic
and mixed characteristic) there exists some $\hat  r \geq 1 $ such
that $p^{\hat r}$ kills $\t{H}^1(\Gamma, M')$. On the other hand,
note that  $\#(\Gamma) $ (the order of $\Gamma$) also kills
$\t{H}^1(\Gamma, M')$. Thus, in the mixed characteristic case,
$p^{r_0}\in(\pi)^r R^\times $ kills $\t{H}^1(\Gamma, M')$.

Now let us deal with  the equi-characteristic case, which needs a
much more complicated argument. We begin by treating the case $M'=
R'$. First let us reduce the problem to the case that $\Gamma$ is
cyclic with order $p$.

Note that if $\#(\Gamma)$ is prime to $p$, since $p$ also kills
$\t{H}^1(\Gamma, R')$  we have $\t{H}^1(\Gamma, R')=0$. On the
other hand, since $\Gamma$ is solvable,  by the
inflation-restriction exact sequence
\begin{equation}\label{inflation}
\begin{split}
\xymatrix{0 \ar[r] & \t{H}^1(\Gamma/\Gamma', {R '}^{\Gamma '} )
\ar[r]^-{\t{Inf}} & \t{H}^1 (\Gamma, R') \ar[r]^-{\t{Res}}
&\t{H}^1(\Gamma', R'),  }
 \end{split}
\end{equation}
where $\Gamma'$ is the  wild inertia subgroup of $\Gamma$, we
reduce the problem to the case that the order of $\Gamma$ is a
power of $p$.

Now we can assume $\Gamma' $ and $\Gamma$ in (\ref{inflation}) are
of
  $p$-power order and we assume that the lemma is proven for cyclic groups of order $p$.
We assume by induction on $\#(\Gamma)$ that it is settled for
smaller $p$-power order and we take $\Gamma ' \triangleleft
\Gamma$ a  nontrivial normal subgroup. Let $K_1$ be the fraction
field of $R'^{\Gamma'}$ and $\pi_1$ a uniformizor $R'^{\Gamma'}$.
By induction,   $\pi^{v( \Delta_{K_1/K}) }$ kills
$\t{H}^1(\Gamma/\Gamma', {R '}^{\Gamma '} )$ and
$\pi_1^{v_{1}(\Delta_{K'/K_1})}$ kills $\t{H}^1(\Gamma', R') $,
where $v_1$ is the valuation on $R'$ by putting $v_1(\pi _1 )=1$.
By the transitivity formula for discriminants, $$\pi^{v(
\Delta_{K_1/K})} \cdot \t{N}_{K_1/K}
(\pi_1^{v_{1}(\Delta_{K'/K_1})})=u \pi^{v(\Delta_{K'/K})}$$ where
$u$ is a unit in $R'$. This product kills $\t{H}^1(\Gamma, R') $.
Thus,  it remains to prove that $\pi^{v(\Delta_{K'/K})}$ kills
$\t{H}^1(\Gamma, R')$ when $\Gamma $ is cyclic order $p$.

Now suppose  that $\Gamma$ is cyclic with order $p$. Let $\varphi$
be a cocycle in $\t{H}^1(\Gamma, R')$ and $\sigma $ a generator of
$\Gamma$. Since $\t{H}^1(\Gamma, K')=0$, we have $\varphi(\sigma )
=\sigma(x/\pi^n) - x/\pi^n $ for an $x\in R'$ and some $n\geq 0$.
Write $x=\sum_{i=0}^{p-1}x_i (\pi')^i$, where $x_i \in R $ for
$i=0, \dots, p-1 $. Then $\varphi(\sigma)= \frac{1}{\pi^n
}\sum_{i=1}^{p-1}x_i(\sigma(\pi')^i -(\pi')^i)$. Note that
$\varphi(\sigma)\in R'$, so we have
\begin{equation}\label{2.2}
\sum_{i=1}^{p-1}x_i(\sigma(\pi')^i -(\pi')^i) \equiv 0 \mod \pi ^n
\end{equation}
Recall that $r = v(\Delta_{K'/K}) $ by definition. Now it suffices
to prove that  if $n >r $ then $x_i \equiv 0 \mod \pi^ {n-r} $ for
all $i=1, \dots , p-1$.

Since $\Gamma$ is wild inertia, there  exists a $\tau \in \Gamma $
with $\tau \not = \t{Id}$ such that
$$\tau(\pi') -\pi ' \equiv 0 \mod (\pi')^2$$

From \eqref{2.2} we have
$$\left(\sum_{i=1}^{p-1}x_i(\frac{\tau(\pi')^i
-{\pi'}^i}{\tau(\pi')-\pi'})\right)(\tau(\pi') -\pi ')\equiv 0
\mod \pi ^ n
$$  Since $v(\tau(\pi') -\pi ')\leq v(\Delta_{K'/K})=r$,
\begin{equation}\label{2.3}
\sum_{i=1}^{p-1}x_i(\frac{\tau(\pi')^i
-{\pi'}^i}{\tau(\pi')-\pi'})\equiv 0 \mod \pi ^ {n-r}.
\end{equation} To prove that $x_i \equiv 0 \mod \pi ^{n-r}$, $i =1,
\dots, p-1$, it suffices (by induction on $n>r$) to prove that
$\pi | x_i $ for all $i=1, \dots, p-1$. First from \eqref{2.3}, we
see that $\pi'| x_1$, so $\pi | x_1$. Suppose $j$ is the minimal
such that $\pi\nmid x_j $, or equivalently $\pi' \nmid x_j$. We
prove that this is impossible by the following claim:

Claim:$$  v(\frac{\tau(\pi')^i
-{\pi'}^i}{\tau(\pi')-\pi'})=\frac{i-1}{p},\  i=1, \dots, p-1.$$

Assuming the claim is true for a moment, we have
$$v\left( \sum_{i=1}^{p-1}x_i(\frac{\tau(\pi')^i -{\pi'}^i}{\tau(\pi')-\pi'})\right)=
v(x_j(\frac{\tau(\pi')^j
-{\pi'}^j}{\tau(\pi')-\pi'}))=\frac{j-1}{p}<1
$$ From \eqref{2.3}, we see this is impossible.

Finally, let us compute $v(\frac{\tau(\pi')^i
-{\pi'}^i}{\tau(\pi')-\pi'})$ to prove the Claim. By the binomial
expansion,

$$ \tau(\pi')^i- {\pi'}^i = (\pi' +(\tau (\pi')-\pi'))^i-{\pi'}^i= \sum_{m=1}^{i} {{i}\choose {m }} (\pi')^{i-m}(\tau(\pi')
-\pi')^m.
$$
 Since $v(\tau(\pi')-\pi')\geq
v(\pi')^2$, we have $v(\frac{\tau(\pi')^i
-{\pi'}^i}{\tau(\pi')-\pi'})= v(i {\pi'}^{i-1})=\frac{i-1}{p}$.
This proves the Claim.

Now let us treat the general case for $M'$ a finite free
$R'$-module rank $m$.  The case $r_0=0$ is clear since $\#(\Gamma)
\in R^\times$ in such case. Now consider the case $r_0>0$. Let us
first reduce the problem to the case that $M'$ is free $R'$-rank
1.
 In fact,  we have Galois descent data on  $M'
\otimes_{R'} K' $ from $\phi_\sigma$.  By ordinary Galois descent,
$M' \otimes_{R'} K'$ has an effective descent. Thus we can find a
nonzero $\alpha \in M' $  such that $\phi_\sigma(\alpha)= \alpha$.
Let $M'_1$ be the saturation of $R'\cdot \alpha \subset M'$;  that
is, $M'_1$ is the smallest submodule of $M'$
 containing  $R'\cdot \alpha$ such that  $M'/M'_1$ is a free $R'$-module.
It is easy to see that  $M'_1$ is a Galois-stable submodule of
$M'$. Thus, we have an exact sequence of a $\Gamma$-descent
modules,
$$0 \to M'_1\to M' \to M'/M'_1 \to 0,   $$
 where $M'/M'_1 $ is a finite free $\Gamma$-descent module with $R'$-rank $m-1$. By  induction on $m$,
 it now suffices to prove the lemma for the case when $M'$ is free
 $R'$-module with rank 1.

Now assume $M'$ is free rank-1 $R'$-module. Select a basis $\alpha
\in M'$. The semi-linear Galois structure on $M'$ is determined by
the 1-cocycle $\sigma \mapsto \phi_\sigma (\alpha)/\alpha \in
R'^\times $. By Hilbert's Theorem 90, we see there exists an $x
\in K'^\times $ such that $\sigma(x)/x =  \phi_\sigma
(\alpha)/\alpha $ for all $\sigma \in \Gamma$.  After multiplying
some power of $\pi$, we can assume $0 \leq v(x)< 1$. Define an
$R'$-module morphism $f_x : M' \to R'\cdot x \subset R'$ by
$\alpha\mapsto x $.  Since $\sigma(x)/x =  \phi_\sigma
(\alpha)/\alpha $, we see that $f_x$ is Galois-equivariant. From
the short exact sequence
$$\xymatrix{0\ar[r]& M' \ar[r]^-{f_x} & R' \ar[r] & R'/
f_x(M')\ar[r]& 0 ,  }$$ we get a short exact sequence $(R'/
f_x(M') )^ \Gamma \to \t{H}^1(\Gamma, M') \to \t{H}^1(\Gamma,
R')$. Note that $f_x(M')=xR'$ and $0 \leq v(x)<1$, so $\pi $ kills
$(R'/ f_x(M') )^ \Gamma$.  As we have seen above, $\pi ^ r$ kills
$\t{H}^1(\Gamma, R')$, so  $\pi^{r+1}$ kills $\t{H}^1(\Gamma,
M')$.

\end{proof}
\begin{remark} Of course, the above proof for the equi-characteristic case also works for
the mixed-characteristic case. Unfortunately, $v(\Delta_{K'/K})$
is always larger than $v(p)r_0$ in this case, so we use $v(p)r_0$
instead of  $v(\Delta_{K'/K})$.  This suggests that there may
still be room to improve the constant $r$.

\end{remark}

\begin{proof}[Proof of Proposition \ref{func}]
 Since $R'/\pi^{n}$ is faithfully flat over
$R/\pi^{n}$ and $M_i'$ is a finite free $R'/\pi^ n $-module,  $M_i
$ is a finite flat (hence free)$R/\pi^{n}$-module.

Now  consider the following commutative diagram:
$$
\xymatrix{
M_1\ar@{^{(}->}[r]^-{i}&(M'_1)^{\Gamma}\ar[d]^{\t{pr}_1^\Gamma}\ar[r]^{f^\Gamma}
&(M_2')^\Gamma\ar[d]^{\t{pr}_2^\Gamma}\\
& (M'_1)^{\Gamma}\ar[r]^{ {f}^\Gamma}& (M'_2)^\Gamma }
$$
where $\rm{pr}_1^\Gamma $,  $\t{pr}_2^\Gamma$ are induced by the
morphism $\pi^{r}$ and $i$ is the natural embedding. Since $M'_i$
is a finite free $R'/\pi^n$-module,
 the morphism $f^\Gamma\circ
\t{pr}_1^\Gamma\circ i $ is an effective descent of $\pi^{r} f $
by the following lemma. The uniqueness of this descent is clear
because $R'/\pi^n$ is flat over $R/\pi^n$.
\end{proof}

\begin{lemma}
Let  $\t{pr}^\Gamma: (R'/\pi^n)^\Gamma \to (R'/\pi^{n})^\Gamma$ be
the morphism  induced by  $\pi^r$. The image of  $\t{pr}^\Gamma$
is $\pi^r R/\pi^{n}$.
\end{lemma}
\begin{proof}  Consider the following commutative diagram of exact sequences:
 $$
\xymatrix{ 0\ar[r] & R'^\Gamma\ar[r]^{\pi^{n}}\ar[d]^{\pi^{{r}}} &
R'^\Gamma\ar[r] \ar[d]^{\t{Id}}  & (R'/\pi^{n}R')^\Gamma
\ar[d]^{h}\ar[r]& \text{H}^1(\Gamma, R')[\pi^{n}]\ar[d]^{\pi^r}\ar[r]&0\\
0\ar[r] & R'^\Gamma \ar[r]^{{\pi^{n-r}}} &R'^\Gamma \ar[r]^-{j} &
(R'/\pi^{n-r}R')^\Gamma \ar[r] &
\text{H}^1(\Gamma, R')[\pi^{n-r}]\ar[r] &0 \\
}
$$
where $\text{H}^1(\Gamma, R')[\pi^{n}]$ is the kernel of $\pi^n$
on $\text{H}^1(\Gamma, R')$ and $h$ is induced by the canonical
reduction map. We need to show $h$ has image $R/\pi^{n-r}=
\t{Im}(j)$. By Lemma \ref{bound},  $\pi^r$ kills
$\text{H}^1(\Gamma, R')$, so
$\t{Im}(\t{pr}^\Gamma)=\t{Im}(j)=(R'^\Gamma)/\pi^{n-r}=R/\pi^{n-r}$.
\end{proof}

\begin{con} \label{bla}
Unfortunately, from now on we cannot deal with the mixed
characteristic case and equi-characteristic case in the same
relative setting. More precisely, in the mixed characteristic case
we can only state and prove results for ``$p^n$-torsion" instead
of ``$\pi^n$-torsion" as in the equi-characteristic case, although
the results and ideas of proofs are almost the same if we replace
``$p$"  by ``$\pi$". To save space,
 we will state the results in mixed characteristic case and use ``(resp. )" to indicate the respective results in
equi-characteristic.
\end{con}

\subsection{Criterion for  effectiveness of descent} Now we will
show that for a  given Galois group, the effectiveness of a finite
free $\Gamma$-descent $R'$-module $(M', \phi_\sigma) $ will be
determined by working at a large torsion level determined by $R'$
and $R$. Recall that $r_0$ denotes  the $p$-index of $\Gamma$ (i.e.,
$p^{r_0}\parallel  \#(\Gamma)$).  Let ${d}=\t{Max}(1, 2r_0)$ in the
mixed characteristic case. In the equi-characteristic case,
  let ${d}=1 $ if $r_0=0$ and ${d}= 2m(r+1) $
if $r_0\geq 1$,  where $m$ is the $R'$-rank of $M'$.
\begin{prop}{\label{descent}}
Let  $(M', \phi_\sigma)$ be a $\Gamma$-descent module with $M'$ a
finite free $R'$-module of rank $m$. Then $(M', \phi_\sigma)$ has
an effective descent if and only if $(M'/p^{ d}, \phi_\sigma)$
\t{(}resp. $(M'/\pi^{ d}, \phi_\sigma)$ \t{)} has an effective
descent to an $R/p^{ d}$-\t{(}resp. $R/\pi^d$-\t{)}module
\t{(}necessarily finite and free\t{)}.
\end{prop}

\begin{proof}The ``only if " direction is trivial. For the
converse, let us  first  prove the proposition in the mixed
characteristic case. First consider  the case  ${ d}=2r_0$. Note
$r_0= r/v(p)$ in this case, so  $p^{r_0}$ kills
$\t{H}^1(\Gamma,M')$ by Lemma \ref{bound}.
  Consider following
diagram:
\begin{equation}\label{e}
\begin{split}
\xymatrix {0 \ar[r] &M' \ar[r]^{p^{r_0}} &M' \ar[r] &{M'/p^{r_0}} \ar[r] &0\\
                   0\ar[r] &M'^{\Gamma}\ar@{^{(}->}[u]\ar[r]^{p^{r_0}} & M'^{\Gamma }
\ar@{^{(}->}[u] \ar[r]^-{\alpha} & ({ M'/p^{r_0}})^{\Gamma
}\ar@{^{(}->}[u]\ar[r]^-{\beta}
& \text{H}^{1} (\Gamma, M') [p^{r_0}] \ar[r] &0\\
}
\end{split}
\end{equation}
where $\text{H}^{1} (\Gamma, M') [p^{r_0}] $ is the kernel of
$p^{r_0}$. Let $M=M'^{\Gamma}$ and let $\t{Im}(\alpha)$ denote the
image of map $\alpha$. By Galois descent theory for $K'/K$, $M$ is
a finite free $R$-module with the same rank as that of $M'$ over
$R'$.
 Obviously, $\alpha(M)\simeq M/p^{r_0}M$.

 We now claim that $\t{Im}(\alpha)$ is an effective descent of ${M'/p^{r_0}}$; that is, the natural $\Gamma$-equivariant
map $ \alpha(M)\otimes_R R'\to {M'/p^{r_0}}$ is an isomorphism.

To prove this  claim, let us consider the  following  commutative
diagram:

$$
\xymatrix{
0\ar[r] & M'\ar[r]^{p^{2{r_0}}}\ar[d]^{p^{r_0}} & M'\ar[r] \ar[d]^{\t{Id}}  & M'/p^{2{r_0}}M' \ar[d]^i\ar[r]& 0\\
0\ar[r] & M'\ar[r]^{{p^{r_0}}}        &M' \ar[r] & M'/p^{{r_0}}M' \ar[r] &0 \\
}
$$

Take Galois invariants to  get a new diagram, with $M$ denoting
$M'^{\Gamma}$:
$$
\xymatrix{ 0\ar[r] & M\ar[r]^{p^{2{r_0}}}\ar[d]^{p^{r_0}} &
M\ar[r] \ar[d]^{\t{Id}}  & (M'/p^{2{r_0}}M')^\Gamma
\ar[d]^{i^\Gamma}\ar[r]&
\text{H}^1(\Gamma, M')[p^{2{r_0}}]\ar[d]^{p^{r_0}}\ar[r]&0\\
0\ar[r] & M\ar[r]^{{p^{r_0}}}        &M \ar[r]^-{\alpha} &
(M'/p^{{r_0}}M')^\Gamma
\ar[r]^-{\beta} & \text{H}^1(\Gamma, M')[p^{r_0}]\ar[r] &0 \\
}
$$

Since  $p^{r_0} $ kills $\text{H}^1(\Gamma, M')$,  the right
square of the above  diagram implies
 $\t{Im}(i^{\Gamma})\subseteq \t{Ker}(\beta) $ and the middle square gives us $\t{Im}(\alpha)\subseteq \t{Im}(i^\Gamma)$.
 By the  exactness of $\alpha$ and
$\beta$, we have  $\t{Im}(\alpha)= \t{Ker} (\beta)$. Thus
$\t{Im}(\alpha)=\t{Im}(i^\Gamma)$.

By hypothesis, since $2{r_0} ={d} $, $M'/p ^{2{r_0} }M' $ has an
effective descent to an $R'/{p ^{2{r_0} }}$-module which is
necessarily  flat. Hence, $M'/p ^{2{r_0} }M' $ is spanned by
$(M'/p ^{2{r_0} }M')^\Gamma $ over $R'$. Thus, $M'/p ^{{r_0} }M' $
is spanned over $R'$ by the image of $(M'/p ^{2{r_0} }M') ^\Gamma
$ in $(M'/p ^{{r_0} }M')^\Gamma $, hence by $\t{Im}(\alpha )
\simeq M /p^{{r_0}}M $. Thus, the natural map $\alpha(M)\otimes_R
R'\to M'/p^{r_0} M'$ is surjective. Both sides are finite free
$R'/p ^{r_0} $-modules with the
 \emph {same} $R'/p^{r_0}$-rank, so this natural map is an isomorphism.
This completes the proof that  ${\alpha}(M )$ is an effective
descent of $M'/p^ {r_0} M'$.

Now let us consider the following diagram
$$
\xymatrix{
0\ar[r] & M'\ar[r]^{p^{r_0}} & M'\ar[r]  & M'/p^{{r_0}}M' \ar[r] & 0\\
0\ar[r] & M\otimes_R R' \ar[u]\ar[r]^{{p^{r_0}}} &M\otimes_R R'
\ar[r]\ar[u]^h &
 (M/p^{{r_0}}M)\otimes_R R'\ar[u]^g \ar[r]  &0 \\
}
$$
where $M=M'^{\Gamma}$ and $h, g$ are the natural maps. We have
proved that the  natural map $g$ is an isomorphism, so by
Nakayama's lemma  $h$ is surjective. Since $M$ is a finite free
$R$-module with  $R$-rank equal to the $R'$-rank of $M'$,  $h$
must be an isomorphism. This settles the proof when $d=2r_0$ in
the mixed characteristic case.

In the mixed characteristic case with ${d}=1$, since
$\t{H}^1(\Gamma, M')=0$ the usual projection argument shows that
 $M/p \simeq (M'/ p)^\Gamma $.
The hypothesis on $M'/ p $ having an effective descent implies
that $(M'/p)^\Gamma \otimes_{R/p} R'/p \to M'/p $ is surjective,
so $R'\otimes_R M \to M' $ induces a surjective map modulo $p$ and
hence is surjective. This surjection  must be an isomorphism for
freeness reason, as above.

Using Lemma \ref{bound}, the proof for the equi-characteristic
case is  the same as above if we replace $p$ by $\pi$ and $r_0$ by
$m(r+1)$ everywhere.
\end{proof}

\section{The Deformation Lemma }
We will  use  results of Grothendieck and Drinfeld to show that
under the hypothesis  of the Theorem \ref{th1},  the effective
descent for $(\G', \w \phi_\sigma)$ does exist at some
$R/\pi^n$-level (see Theorem \ref{deformation} below). First we
need to recall several facts concerning  deformations of
\emph{truncated \bt groups}.

\subsection{Deformations theory from Grothendieck and Drinfeld}
Recall that   a \emph{truncated \bt group of level $N\geq 2$ with
height $h$} over a scheme $S$ is an inductive system $(G_n,
i_n)_{0\leq n\leq N}$ of finite locally free commutative $S$-group
schemes and closed immersions $i_n: G_n\to G_{n+1}$ such that
\begin{enumerate}
\item $|G_n|=p^{nh}$ \item $i_n $ identifies $G_n$ with
$G_{n+1}[p^n]$ for $0\leq n<N$.
\end{enumerate}
\emph{A truncated \bt group of level 1 with height $h$} is finite
locally free $S$-group scheme killed by $p$ with order $p^h$ such
that the sequence
$$\xymatrix{G_0 \ar[r]^{F}& G_0^{(p)}\ar[r]^{V}& G_0}$$
is exact over $\t{Spec}(\O_S/p\O_S)$, where $F$ is the { relative}
``Frobenius" and $V$ is ``Verschiebung" and $G_0 = G\times_S
\t{Spec}(\O_S/p\O_S)$.


Now let us recall a result of Grothendieck on truncated \bt
groups.

\begin{theorem}[Grothendieck]\label{Grothendieck}
Let $i: S\to S' $ be a nilimmersion of schemes and $G$ a truncated
\bt group over $S$ of level $n\geq 1$. Suppose $S'$ is affine.
\begin{enumerate}
\item There exists a truncated \bt group $G_n'$ over $S'$ of level
$n$ such that $G'_n$ is  a deformation of $G_n$.
 \item If
$S$ is complete local noetherian ring  with perfect residue field
of  characteristic $p$, then there exists a \bt group $H$
 such that $G=H[p^n]$.
 \item If there exists a \bt  group $H$ over $S$ such that $G=H[p^n]$, then for any deformation of $G'$ of $G$ over
$S'$, there exists a \bt group $H'$ over $S'$ such that $H'$ is a
deformation of $H$ and $G'\simeq H[p^n]$ \t{(}as a deformation of
$G$\t{)}.

\end{enumerate}

\end{theorem}

\begin{proof}See pp. 171-178, \cite{luc}.
\end{proof}

We next review  Drinfeld's results (\S1, \cite{katz}) on
deformations of \bt groups.  Let $A$ be a ring, $N\geq 1 $ an
integer such that $N$ kills $A$, and $I$  an ideal in $ A$
satisfying  $ I^{\mu+1}=0$ for a fixed $\mu >1$. Let $ A_0 =A/I $.
\begin{lemma}[Drinfeld]\label{drinfeld}
Let $G$ and $H$ be  \bt groups over $A$.
\begin{enumerate}
\item The groups $\textnormal {Hom}_{A\t{-gp}}(G, H)$ and
$\textnormal {Hom}_{A_0\textnormal{-gp}}(G_0, H_0)$ have no
non-trivial $N$-torsion. \item The natural map  $ \textnormal
{Hom}_{A\t{-gp}}(G, H)\rightarrow \textnormal{Hom}_{A_0\t{-gp}}
(G_0, H_0)$ is injective. \item For any homomorphism $f_0:
G_0\rightarrow H_0 $, there is a unique homomorphism $``N^{\mu}f"
: G\rightarrow H $
  which lifts $N^{\mu}f_0$.
\item  In order that a homomorphism $f_0 : G_0 \to H_0 $ lift to a
\t{(}necessarily unique\t{)} homomorphism $f: G\to H$, it is
necessary and sufficient that  the homomorphism $``N^{\mu}f" :
G\to  H$ annihilate
 the subgroup $G[N^{\mu}]$
of $G$.
\end{enumerate}
\end{lemma}
\begin{proof} See \S1, \cite{katz}.\end{proof}

\begin{remark} Drinfeld actually shows the above results  for more general \emph{fppf group sheaves},
including  \bt groups as a special case.
 \end{remark}
Now let us  refine Lemma \ref{drinfeld} for our purposes.

\begin{lemma} {\label{defor}}
Keep the notations in Lemma \t{\ref{drinfeld}} with  $N=p^m$ and
$\mu=1$.  Choose an integer $l\geq 1$. Let $f_0: G_0\to H_0$ be a
morphism of \bt groups over $A_0$. If
 there
exists an $A$-group morphism $g:G[p^{ m+l }]\to H[p^{ m+l} ]$ such
that $g\otimes_A A_0= f_0|_{G_0[p^{m+l}]}$, then there exists an
$A$-morphism $f: G\to H $ such that $f\otimes_A A_0 =f_0$ and
$f|_{G[p^l]}=g|_{G[p^l]}$. Furthermore, if $g $ and $f_0$ are
isomorphisms, then $f$ is an isomorphism.
\end{lemma}

\begin{proof}  Using  Lemma \ref{drinfeld} in our case, we first prove we can
 lift the morphism $f_0$ to an $A$-morphism. By Lemma {\ref{drinfeld}} (4), it is
necessary that the lift of $``p^m f_0"$
 annihilates  the subgroup $G[p^m]$.
 Let  $B$ be any $A$-algebra, $B_0=B/IB$ and $F=``p^mf_0"$.   By \S1, \cite{katz},
 $F$ is defined in the following (using formal smoothness of $H$):
\begin{equation}\label{F1}
\begin{split}
\xymatrix{ G(B) \ar[dr]\ar@{--}[r]&- \ar@{--}[r]^{F}&-\ar@{-->}[r]  & H(B)\\
  &G_0(B_0) \ar[r]^{f_{0}} &H_0(B_0) \ar[ur]_{p^m \times  \t{any lifting}} &
}
\end{split}
\end{equation}
Since $g\otimes_A A_0= f_0|_{G_0[p^{m+l}]}$, we have the following
commutative diagram:
\begin{equation}\label{F2}
\begin{split}
\xymatrix{ G[p^{m+l}](B) \ar@{->}[d]\ar[r]^{g}
&H[p^{m+l}](B)\ar@{->}[d]\\
G_{0}[p^{m+l}](B_0)\ar[r]^{f_{0}} & H_{0}[p^{m+l}](B_0) }
\end{split}
\end{equation}
where the columns are reduction modulo $I$.   Now for all $ x \in
G[p^{m} ](B)$, $g(x)$ is a lift of  $f_{0} ( x \mod I )$. Note
that $g(x)\in H[p^m](B )$, so $F(x)=p^m(g(x))=0$.  Thus, $F=``p^m
f_{0}"$ kills $G[p^m]$, so we can define $f: G\simeq G/G[p^m]\to
H$ induced by $F$ and thus  ${f: G\to H}$  lifts $f_0$ and $p^mf
=F$. Tracing the definition of $F$ in (\ref{F1}) and
 (\ref{F2}), we get $F|_{G[p^{m+l}]}=p^m g$. Thus $f, g : G[p^l ]\rightrightarrows H[p^l]$
 have the same composite with the {\it fppf} map $p^m : G[p^{m+l}] \to G[p^l]$, so
 $f|_{G[p^l]}=g|_{G[p^l]}$.

If $f_0$ and $g$ are isomorphisms, we can use $f_0^{-1}$ and
$g^{-1}$ to construct an $A$-morphism of \bt groups $\tilde f: H
\to G $  which lifts $f_0^{-1}$. We see that $(f\circ \tilde
f)\otimes_A A_0$ and $(\tilde f \circ f ) \otimes_A A_0 $ are
identity maps.  By Lemma \ref{drinfeld} (4), $f\circ \tilde f $
and $ \tilde f \circ f  $ have to be identity maps,  so  $f$ is an
isomorphism.

 \end{proof}

\subsection{Application to descent}

 We use notations as in Theorem \ref{th1}, so  $(\G', \w \phi_\sigma)$ is an $R'$-\bt group
equipped with $\Gamma$-descent datum
 constructed from $G\otimes_K {K'}$ over $K'$ as in Introduction. As in  Convention \ref{bla},
we will state  results in the  mixed characteristic case and use
``(resp. )" to indicate the corresponding results in the
equi-characteristic case.
 Let $\lambda=1$ for $p\not =2$ and $\lambda=2$ for $p=2$.  For any real number $x$, define
$\left  <x\right> $ the maximal integer $< x$,  so $\left <
x\right>= x-1$ if $x\in \Z$.  Put $\nu= \left <\log_2 e\right
> +\lambda$ in the mixed characteristic case and $\nu= \left
<\log_2 \hat e\right > +\lambda $  in the equi-characteristic
case,
 where $\hat e$ is  the relative ramification index $e(K'/K)$.
For each $n\geq 1$,  let  $\mathcal  G '_{n}=\mathcal  G
'\otimes_{R'}R'/p^{n}$ and  $l_n= 3n +\nu $ in the mixed
characteristic case,  and $\mathcal  G '_{n}=\mathcal  G
'\otimes_{R'}R'/\pi^{n}$ and
 $l_n= \left <\log_2n\right >+1 +\nu  $ in the equi-characteristic case.
    We will use Grothendieck's theorem and the same argument of Drinfeld to prove the following:

\begin{theorem}[Deformation Lemma]\label{deformation}
If   $\mathcal G'[p^{ l_n}]$ has an effective descent to a
truncated \bt group $\mathcal H $ over $R$, then there exists a
\bt group  $\mathcal G_{n}$ over $R/p^n$ \t{(}resp. $R/\pi^n$\t{)}
 such that $\mathcal  G_{n}$ is an effective descent for  $ \mathcal  G '_{n}$.
\end{theorem}
\begin{remark}
We do not claim to identify with $\H $ with $\G_n[p^{l_n} ] $ as
descents of $\G'_n[p^{l_n}]$.
\end{remark}
\begin{proof} Let $\G'_0=\G'\otimes_{R'}k  $ be the closed fiber
of $\G'$ over $k$.  Since $R$ and $R'$ have  the same residue
field, it is obvious that the closed fibre $\mathcal H_0 =\mathcal
H\otimes_{R} k$ is naturally isomorphic to $\mathcal G'_0[p^{l_n
}]$.  By Theorem \ref{Grothendieck},  there  exists a \bt  group
$\mathcal G_n $ over $R/p^n$ (resp. $R/\pi^n$) such that $\G_n$ is
a deformation of $\G'_0$ over $R/p^n$ (resp. $R/\pi^n$) and $\H_n
=\G_n[p^{l_n}]$ as a deformation of $\mathcal H_0 \simeq
\G_0'[p^{l_n}]$, where $\H_n = \H \otimes_R {R/p^n}$ (resp. $\H_n
= \H \otimes_R {R/\pi ^n}$). It now suffices to prove:

\begin{lemma}\label{tt}
 There exists an $\Gamma$-equivariant isomorphism
 ${f:  \mathcal G'_n  \ito  \mathcal G_n\otimes_R R' }$ as \bt groups  over $R'$ with   $f$
 lifting the identification of
  $\G'_0 $  with the  closed fiber $\G_0$ of  $\G_n$.
\end{lemma}


 For the proof of Lemma \ref{tt}, first note that both $\G'_n$ and $\G_n\otimes_R R' $ are
deformations of $\G'_0$ over $R'/p^n $ (resp. $R'/\pi^n $). Since
$\H_n$ is an effective descent of $\G'_n[ p^ {l_n}]$ and $\H_n
=\G_n[p^{l_n }]$ as deformations of $\H_0 \simeq \G_0 [p ^{l_n
}]$, there exists a Galois-equivariant isomorphism $g: \G'_n[ p^
{l_n}]\ito (\G_n\otimes_R R') [p^{l_n }]$ such that
the isomorphism $f_0 : \G'_0 \ito \G_n\otimes_R k $ agrees with
$g$ on $p^{l_n}$-torsion over $k$.

Now let us first discuss the proof for the case that $R$ has mixed
characteristic. Recall that $l_n= 3n+ \nu=3n+ {\left <\log _2
e\right >} + \lambda $ in this case.   Consider the following
 successive square-zero thickenings, where  $\m'$ is the maximal ideal
of $R'$:
$$k=R'/\m'\leftarrow  R'/{\m'}^2\leftarrow \cdots \leftarrow R'/ {\m'}^{2^s} \leftarrow R'/{\m^\prime}^{2^{s+1}}\leftarrow
\cdots \leftarrow R'/{\m'}^{2^{\left <\log _2 e\right >}}
\leftarrow  R'/pR' $$ By Lemma \ref{defor} and induction on $1\leq
s \leq \left <\log_2 e\right > $ with $N = p $, there exists an
isomorphism $f_1: \G'_1\ito \G_1\otimes_{R/p}{R'/p} $ such that
$f_1|_{\G'_1[p^{3n+\lambda-1}]}= (g\otimes_{R'/p ^n}R'/p
)|_{\G'_1[p^{3n+\lambda-1}]}$.

Consider deformations  of $\G'_1=\G'_n\otimes_{R'/p^n}R'/p$ along
the following:
$$R'/p\leftarrow R'/p^2\leftarrow\cdots\leftarrow R'/p^{2^s}\leftarrow R'/p^{2^{s+1}}\leftarrow \cdots\leftarrow
R'/p^{2^{\left <\log_2n \right >}} \leftarrow R'/p^n.$$ By
induction on $1\leq s \leq \left <\log _2n \right >$ and Lemma
\ref{defor} with $N=p^{2^s }$ ($N=p ^n$ in the last step), there
exists an isomorphism $f_n: \G'_n \ito \G_n \otimes_{R/p^n}R'/p^n
$ such that $f_n|_{\G'_n[p^{3n-1+\lambda-t}]}=
g|_{\G'_n[p^{3n-1+\lambda-t}]}$, where $t=(\sum_{i=1}^{\left
<\log_2{n}\right >}2^i)+n = 2 ^{\left < \log_2 n\right>+1} -2+n
\leq 3n-1 $. Thus we have $f_n|_{\G'_n[p^{\lambda}]}=
g|_{\G'_n[p^{\lambda}]}$.

In the equi-characteristic case,  we will also consider
deformations of $\G_0 $ along
$$k=R'/\m'\leftarrow  R'/{\m'}^2\leftarrow \cdots \leftarrow R'/ {\m'}^{2^s} \leftarrow R'/{\m^\prime}^{2^{s+1}}\leftarrow
\cdots \leftarrow R'/{\m'}^{2^{\left <\log _2 e\right >}}
\leftarrow  R'/\pi R' $$
 and  deformations of $\G'_1= \G'_n \otimes _{R'/ \pi^n }R'/\pi$ along
$$R'/\pi\leftarrow R'/\pi^2\leftarrow\cdots\leftarrow R'/\pi^{2^s}\leftarrow R'/\pi^{2^{s+1}}\leftarrow \cdots\leftarrow
R'/\pi^{2^{\left <\log_2n \right >}} \leftarrow R'/\pi^n. $$ We
also use Lemma \ref{defor} to do  induction. The only difference
from the mixed characteristic case is that $N$ is always $p$.
Similarly, we get an isomorphism $f_n: \G'_n \ito \G_n
\otimes_{R/\pi^n}R'/\pi^n $ such that  $
f_n|_{\G'_n[p^{\lambda}]}= g|_{\G'_n[p^{\lambda}]}$.

In both the mixed characteristic and equi-characteristic cases,
let $f=f_n$. It suffices to prove that the $R'$-isomorphism of \bt
groups ${f: {\mathcal G'_n}\ito {\mathcal G_n}\otimes_R R' } $  is
Galois equivariant. By Lemma \ref{check} and  Lemma \ref{drinfeld}
(2), it suffices to check $f_0 =(f_0)_\sigma$ on $k$-fibers. Note
that $\psi : {\sigma \mapsto  f_0^{-1}\circ ( f_0)_\sigma }$
define a cocycle in ${C}^1(\Gamma, \t{Aut}_ {R'}{(\G'_0)})$. But
since $R'$ is totally ramified over $R$, the $\Gamma$-action on
$k'=k$ is trivial.  Thus,  $\psi $   induces a representation
$\Gamma \to  \text{Aut}_k(\mathcal G'_0) $. Now from the condition
that $f|_{\G'_n[p^{\lambda}]}= g|_{\G'_n[p^{\lambda}]}$ with  $g:
\G'_n[p^{l_n} ] \to \G_n[p ^{l_n }]\otimes_R R '  $ a
$\Gamma$-equivariant morphism using the $\Gamma$-descent datum for
$\G'_n[p^{l_n}] $, we see that $f|_{\G'_0[p^{\lambda}]}$ is
Galois-equivariant. Thus, via  $\psi $ the action by $\Gamma$ on
$\G'_0[p^{\lambda}]$ is trivial. Now we can conclude that $ \psi$
is trivial via the following lemma.
\end{proof}
\begin{lemma}
 Let $\psi \in \textnormal{Aut}_k(\mathcal G'_0 )$ be an automorphism of finite order. If
$\psi|_{\mathcal G'_0[p^\lambda]}$ is trivial \t{(}recall
$\lambda=1$ for $p>2 $ and $\lambda =2$ for $p=2$\t{)}, then
$\psi$ is trivial.
\end{lemma}
\begin{proof} Passing to the finite order $W(k)$-linear automorphism  $\mathbb D(\psi)$ of the Dieudonn\'e  module
$ \mathbb D(\G_0)$ as a finite free $W(k)$-automorphism, the
convergence property of $p$-adic logarithm on
$\t{GL}_{W(k)}(\mathbb D(\G_0))$ gives the result.
\end{proof}

\section{Proof of the Descent Theorem }
\subsection{Messing's theory on deformations of \bt groups}
 We now prepare to prove Theorem \ref{th1}. Let us first recall some results from \cite{mess} on the
 deformation theory of \bt groups.

Let $S$ be an affine scheme on which  $p$ is nilpotent,  so  there
exists an integer $N>0$ such that $p^N$ kills $\O_S$. Let $\G$ be
a \bt group over $S$.

As an {\it fppf} abelian sheaf,  there exists an extension
\begin{equation}\label{first}
{0\to V(\G) \to E(\G) \to \G\to  0 }
\end{equation}
of $\G$ by a vector group $V(\G)$ (i.e., the group scheme
represented by a quasi-coherent locally  free sheaf of finite
rank). The extension is universal in the sense that for  any
extension $0 \to M\to \mathcal E  \to \G \to 0$ of $\G$ by another
vector group $M$ there is a unique linear map $\xymatrix{V(\G)
\ar[r]^-{\phi}& M } $ such that the extension is (uniquely)
isomorphic to the pushout of (\ref{first}) by $\phi$.

In general,  we can define a ``Lie" functor from the category of
group sheaves  over $S$ to the category of vector group schemes
over $S$ (see \S2, Ch. 3, \cite{mess}). Taking ``Lie" of the
universal extension (\ref{first}), Proposition 1.22 in \cite{mess}
ensures that the following sequence of vector groups over $S$ is
exact:
$$ {0\to V(\G) \to \t{Lie}({E}(\G))\to \t{Lie}(\G)\to  0 }.  $$

For a scheme $X$, recall that its crystalline site {\bf Crys}($X$)
consists of the category whose objects are pairs $(U\inj T,
\gamma) $ where:
\begin{enumerate}
\item $U$ is an open sub-scheme of $X$ \item $U \inj T$ is a
locally nilpotent immersion \item $\gamma = (\gamma_n) $ are
divided powers  (locally nilpotent) on the ideal $I$ of $U$ in $T$
\end{enumerate}
A morphism from $(U \inj T , \gamma )$ to $(U' \inj T', \delta )$
is  a commutative diagram
$$\xymatrix{U\ar[d]^f  \ar@{^{(}->}[r] & T \ar[d]^{\bar f }\\
U' \ar@{^{(}->}[r] & T' }
$$
such that $U \to U'$ is the inclusion and $\bar f : T \to T'$ is a
morphism compatible with divided powers.
\begin{remark}
The crystalline site defined above is weaker that  what Berthelot
defined, but it suffices to study the deformation theory of \bt
groups in \cite{mess}. We always drop the notation of divided
powers $ \gamma$ from the pair   $(U\inj T, \gamma) $ if no
confusion will arise.
\end{remark}
Let $I$ be  a quasi-coherent ideal of $\mathcal O_S$ endowed with
locally nilpotent divided powers. Let $S_0 $ be Spec$(\mathcal O_S
/I)$, so $(S_0 \hookrightarrow S)$ is an object of the crystalline
site of $S_0$. Let $\G_0$, $\H_0$ be  \bt groups over $S_0$ and
$\G$, $\H$ respective   deformations of $\G_0$, $\H_0$ over $S$.
(Such deformations always exist  when $S_0$ is affine,  by
Grothendieck's Theorem (Theorem \ref{Grothendieck}), because  $S$
is necessarily affine when $S_0$ is affine.) Let $u_0 : \G_0 \to
\H_0$ be a morphism of \bt groups. By the construction of
universal extension (\ref{first}), there exists natural morphism
$E(u_0): E(\G _0 ) \to E(\H_0)$ induced by $u_0$. The following
results are the technical heart of \cite{mess}:
\begin{lemma}\label{4.2}Keep  notations as above.
\begin{enumerate}
\item There exists a unique morphism $E(u): E(G)\to E(H)$ which
lifts $E(u_0)$. \item Let $\K$ be a third \bt group on $S$ and
$u'_0:  \H_0 \to \K_0$ a morphism. Then $E(u'_0 \circ u_0 )= E(u_0
') \circ E(u_0)$.
\end{enumerate}
\end{lemma}
\begin{proof}See Chap. IV, \S 2, \cite{mess}.
\end{proof}
\begin{co}\label{4.2.5}
Keep notations as above.
\begin{enumerate}
\item If $\G=\H$ and $u_0=\t{Id}_{\G_0}$, then $E(u)=\t{Id}_{\G}$
\item If $u_0$ is an isomorphism, so is  $E(u)$.
\end{enumerate}
\end{co}

For $ \G _0$ a \bt group over $S_0$, let $\D(\G_0)$ denote
$\t{Lie}(E(\G_0))$.  Assuming for a moment that there exists a
deformation $\G$ of $\G_0$ over $S$, let $\D(\G_0)_S$  denote
$\t{Lie}( {E(\G)})$.  By Lemma \ref{4.2}, the $\O_S$-module  $
\D(\G_0)_S $ does not depend on the choice of deformation of
$\G_0$ in the sense that if $\G_1 $ and $\G_2 $ are two
deformations of $\G_0$ over $S$ and $ E(u): E(\G_1) \to E(\G_2)$
is the unique isomorphism lifting the identification of $E(\G_0)$
over $S_0$, then $\t{Lie}(E(u)) $ is a natural isomorphism from
$\t{Lie}(E(\G_1))$ to $\t{Lie} (E(\G_2))$. In this way, we  glue
the vector bundles  $\t{Lie}(E(\G))$ over open affines in $S$
(where deformations of $\G_0$ do exist!) to construct a crystal
$S\mapsto \mathbb \D(\G_0)_S$ on $S_0$ even in the absence of
global lifts.

An  $\O_S$-submodule  $\t{Fil}^1 \hookrightarrow \D(\G_0)_S $ is
said to be \emph{admissible} if $\t{Fil}^1 $ is a locally-free
vector subgroup scheme over $S$ with locally-free quotient such
that it  reduces to $V(\G_0) \hookrightarrow \t{Lie}(E(\G_0)) $
over  $S_0$.

Define a category $\mathscr C$ whose objects are pairs $(\G_0,
\t{Fil}^1)$ with $\G_0 $ a \bt group over $S_0$ and $\t{Fil}^1$ an
admissible filtration on $\D(\G_0)_S$. Morphisms are defined as
pairs $(u_0, \xi)$ where $u_0: \G_0 \to \G_0'$ morphism of \bt
groups over $S_0$ and $\xi$ is morphism of $\O_S$-modules such
that the diagram
$$
\xymatrix{\t{Fil}^1 \ar@{^{(}->}[r]\ar[d]_{\xi} & \D(\G_0)_S\ar[d]^{\D(u_0)_S}\\
\t{Fil}'^1\ar@{^{(}->}[r] & \D(\G'_0)_S }
$$
commutes and over $S_0$ reduces  to
$$
\xymatrix{V(\G_0) \ar@{^{(}->}[r]\ar[d]_{V(u_0)} & \t{Lie}(E(\G_0))\ar[d]^{\t{Lie}(E(u_0))}\\
V(\G'_0)\ar@{^{(}->}[r] & \t{Lie}(E(\G'_0)) }
$$

\begin{theorem}[Messing]{\label{messing}}
The functor $\G\mapsto (\G_0, V(\G)\hookrightarrow\t{Lie}(E(\G)))$
establishes an equivalence relation of categories:
\begin{center}
\t {\bt groups over }$S$  $\longrightarrow \mathscr C. $
\end{center}
\end{theorem}
\begin{proof} Chap. V, \S 1,  \cite{mess}.  \end{proof}

\subsection{Proof of Theorem \ref{th1} }

Now we shall  use the above theory to prove Theorem \ref{th1}.
Let us first construct the constant $c_0$. Recall  the constants
$d$ and $l_n $ constructed  in Proposition \ref{descent} and
Theorem \ref{deformation} respectively. Set $c_0=l_{d}$ in both
the mixed characteristic case and the equi-characteristic case. As
in Convention \ref{bla}, we will state the proof for the mixed
characteristic case, and use ``(resp. )" to refer the
corresponding result in the equi-characteristic case.

  Consider the following diagram of deformations:
$$
 \xymatrix{
     k=R'/{\m}'\ar@{=}[d]  &  R'/p\ar[l] & R'/p^2 \ar[l]  & \ar[l] \cdots& R'/p^n\ar[l] & R'/p^{n+1}\cdots\ar[l]\\
                  R/\m     &  R/p\ar[u] \ar[l]&  R/p^2\ar[u]\ar[l] &\ar[l]  \cdots &R/p^n\ar[l]\ar[u] &R/p^{n+1}\cdots\ar[l]
\ar[u]\\
            }
$$
Respectively, in the equi-characteristic case,  we consider the
following diagram:
 $$
 \xymatrix{
     k=R'/{\m}'\ar@{=}[d]  &  R'/\pi\ar[l] & R'/\pi^2 \ar[l]  & \ar[l] \cdots& R'/\pi^n\ar[l] & R'/\pi^{n+1}\cdots\ar[l]\\
                  R/\m     &  R/\pi\ar[u] \ar[l]&  R/\pi^2\ar[u]\ar[l] &\ar[l]  \cdots &R/\pi^n\ar[l]\ar[u] &R/\pi^{n+1}\cdots\ar[l]
\ar[u]\\
            }
$$
For each $n\geq 1 $,  let $\mathcal  G'_n =\mathcal  G'\otimes_
{R'} R'/p^n $ (resp. $\mathcal  G'_n =\mathcal  G'\otimes_ {R'}
R'/\pi^n)$.   We have
  universal extensions on each level $n$:
$$0\rightarrow V(\mathcal  G'_n)\rightarrow E(\mathcal  G'_n)\rightarrow \mathcal  G'_n \rightarrow 0.$$

Using the ``Lie" functor, we  get an exact sequence of finite free
$R'/p^n $-modules (resp. $R'/\pi^n$-modules)
$$ 0\rightarrow V(\mathcal  G'_n)\rightarrow \t{Lie} (E(\mathcal  G'_n))\rightarrow \t{Lie} (\mathcal  G'_n) \rightarrow 0.$$

Put  $N'_n=V(\mathcal  G'_n),\  M'_n= \t{Lie} (E(\mathcal G'_n))$,
and $ L'_n= \t{Lie} (\mathcal  G'_n) $. For each $n$, we have an
exact sequence
\begin{equation}\tag{$S'_n$}\label{s'n}
0 \rightarrow N'_n \rightarrow M'_n \rightarrow L'_n\rightarrow 0.
\end{equation}
 By the functorial construction of the universal extension (\ref{first}),
  we see that $N'_n$, $M'_n$ and $L'_n $ are imbued with $\Gamma$-descent data coming from
the generic fiber of $\G'$ for all $n$. Furthermore, ($S'_n$) is a
short exact sequence of $\Gamma$-descent modules and its formation
is compatible with the base changes $ R'/{p^{n+1}}\to R'/p^n$
(resp. $ R'/{\pi^{n+1}}\to R'/\pi^n$).

\begin{lemma}\label{4.3}
 Under the hypothesis of Theorem \t{\ref{th1}}, for each $n\geq 1$
the exact sequence of $\Gamma$-descent modules $(S'_n)$ has an
effective descent
\begin{equation} \tag{$S_n$} \label{sn}
{0\to N_n\to M_n\to L_n\to 0 }
\end{equation}
such that $(S_{n+1})\otimes_{R/ p^{n+1}}R/p^n\simeq (S_n )$
descending $(S'_{n+1})\otimes_{R'/ p^{n+1}}R'/p^n\simeq (S'_n )$
\t{(}resp. $(S_{n+1})\otimes_{R/ \pi^{n+1}}R/\pi^n\simeq (S_n )$
descending $(S'_{n+1})\otimes_{R'/ \pi^{n+1}}R'/\pi^n \simeq (S'_n
)$\t{)}.
\end{lemma}
\begin{proof}
Let  $ N'=  \varprojlim_{n}   N'_n $,  $ M'=  \varprojlim_{n}
M'_n $ and  $ L'=  \varprojlim_{n}   L'_n $. We get a short exact
sequence of $\Gamma$-descent modules:
\begin{equation}\tag{$S'$}\label{s'}
0\to N'\to M'\to L'\to 0,
\end{equation}
 where $N' $, $M'$ and $L'$ are finite free $R'$-modules,  and this  induces each \eqref{sn} by reduction.   It suffices to prove that
 \eqref{s'} has an effective descent.

By the hypothesis of Theorem \ref{th1}, $\mathcal G'[p^{c_0}]$ has
an effective descent $\mathcal H $.  Since $c_0=l_d$, by  Theorem
\ref{deformation} with $n=d$ there exists a \bt group $\mathcal
G_d $ over $R/p^d$ (resp. $R/\pi^d$) such that $\mathcal G_d $  is
an effective descent of $\mathcal G'_d$. If we take the ``Lie"
functor of the crystals constructed for $\mathcal G'_d$ and
$\mathcal G_d $,   we get an exact sequence of $R/p^d$-modules
(resp. $R/\pi^d$-modules)
\begin{equation}\tag{$\tilde S_d $}
{0\to\tilde  N_d\to \tilde M_d \to \tilde L_d\to 0 }
\end{equation}
which is an effective descent of
\begin{equation}\tag{$S'_d$}
0\to N'_d \to M'_d \to L'_d\to 0.
\end{equation}
By Proposition {\ref{descent}}, we see that  $L'$,  $M'$,  and
$N'$ have effective descents to finite free $R$-modules $L$, $M$
and $N$ respectively.  By Corollary \ref{exact},
 $(S')$ has an effective descent to an  exact sequence of finite free
 $R$-modules
\begin{equation}\tag{$S$}
0 \to N\to M\to L\to 0.
\end{equation}
\end{proof}
\begin{remark} \label{remark}We have to be very careful concerning
the following technical point.  As two effective descents of
$(S'_d)$, the sequences $(\tilde S_d)$ and $(S_d)$ are not known
to be isomorphic in the sense that there exists  an $R$-linear
isomorphism of {\it exact sequences}  $f: (\tilde S_d)\to (S_d)$
such that the following diagram commutes $\Gamma$-equivariantly
$$\xymatrix{(\tilde S_d)\otimes_R R'\ar[d]_{f\otimes {\t{Id}}}\ar[r]^-{\tilde i} &(S'_d) \\
(S_d)\otimes_R R' \ar[ur]^{i}& }
$$
where $\tilde i$ and $i$ are isomorphisms of exact sequences of
$\Gamma$-descent modules.  In fact, Example \ref{example} shows
that $(\tilde S_d)$ and $(S_d)$ are not necessarily  isomorphic.
Fortunately,  by Proposition \ref{func}, we see that as effective
descents of $(S'_{d-r_0})$ (resp. $(S'_{d-r})$),
$$ (\tilde S_d)\otimes_R R/p^{d-r_0}\simeq (S_d)\otimes_R R/p^{d-r_0}.  $$
(resp. $(\tilde S_d)\otimes_R R/\pi^{d-r}\simeq (S_d)\otimes_R
R/\pi^{d-r}$).
\end{remark}
Set $m=d-r_0 \geq 1 $ (resp. $m=d-r\geq 1$). We shall  consider
deformations    of $\G_m=\G_d \otimes_{R/p^d}R/p^{d-r_0}$ (resp.
$\G_m=\G_d \otimes_{R/\pi^d}R/\pi^{d-r}$) along
$$\xymatrix{R/p^m &R/p^{2m }\ar[l] &\cdots \ar[l] & R/p^{sm}\ar[l]&R/p^{(s+1)m}\ar[l]&\cdots \ar[l]}$$
(resp. $\xymatrix{R/\pi^m &R/\pi^{2m }\ar[l] &\cdots \ar[l] &
R/\pi^{sm}\ar[l]&R/\pi^{(s+1)m}\ar[l]&\cdots \ar[l]}$)
\begin{lemma}\label{4.4}
 For each $s\geq 1$, $\mathcal G'_{sm}= \mathcal G'\otimes_{R'}
R'/p^{sm}$ \t{(}resp. $\mathcal G'_{sm}= \mathcal G'\otimes_{R'}
R'/\pi^{sm}$\t{)} has an effective descent to an $R/p^{sm}$-
\t{(}resp. $R/\pi^{sm}$-\t{)} \bt group $\mathcal G_{sm}$ and
there are unique isomorphisms $\G_{(s+1)m}\otimes_R R/p^{sm}\simeq
\G_{sm}$ \t{(}resp. $\G_{(s+1)m}\otimes_R R/\pi^{sm}\simeq
\G_{sm}$\t{)} that  descend the isomorphisms
$\G'_{(s+1)m}\otimes_{R'} R'/p^{sm}\simeq \G'_{sm}$ \t{(}resp.
$\G'_{(s+1)m}\otimes_{R'} R'/\pi^{sm}\simeq \G'_{sm}$\t{)}.
\end{lemma}

\begin{proof}
We only prove the mixed characteristic case. In  the
equi-characteristic case the proof  works if we replace $p$ with
$\pi$ everywhere.

 Now let us prove the lemma with by induction on $s$. For $s=1$,
 recall that in the proof of Lemma \ref{4.3}, by the hypothesis of Theorem
 \ref{th1} and Theorem \ref{deformation} with $n=d$  we get an $R/p^d$-\bt group $\G_d $
 which is an  effective
descent of $\G'_d$. Letting $\G_m = \G_d \otimes_{R/p^d}R/p^m$
settles the case $s=1$.
  For  $s=2$, considering   the deformation problem for $\G_m $ along $R/p^m \to R/p^{2m}$,
we need to construct a deformation $\G_{2m}$ of $\G_m$ to
$R/p^{2m}$  such that $\G_{2m}$ satisfies the conditions of Lemma
\ref{4.4}. We claim that it  is equivalent to construct an
admissible filtration
 ${\t{Fil}^1\hookrightarrow
\D(\G_m)_{R/p^{2m}}}$ such that the following diagram commutes
\begin{equation}{\label{diag4.1}}
\begin{split}
\xymatrix{
\t{Fil}^1\otimes_R R'\ar[d]^{f^1}\ar@{^{(}->}[r]& \D(\G_m)_{R/p^{2m}}\otimes_R R'\ar[d]^f\\
V(\G'_{2m})\ar@{^{(}->}[r] & \t{Lie}(E(\G'_{2m})) }
\end{split}
\end{equation}
where $f^1$, $f$ are Galois-equivariant isomorphisms and the
diagram \eqref{diag4.1} lifts the following commutative diagram
\begin{equation}{\label{diag4.2}}
\begin{split}
\xymatrix{
V(\G_m )\otimes_R R'\ar[d]^{f^1_0}\ar@{^{(}->}[r]& \t{Lie}(E(\G_m))\otimes_R R'\ar[d]^{f_0}\\
V(\G'_{m})\ar@{^{(}->}[r] & \t{Lie}(E(\G'_{m})) }
\end{split}
\end{equation}

To make notations easier, we denote the above  diagrams (e.g.,
diagram \eqref{diag4.1}) in the following way:
$$ \xymatrix{(\t{Fil}^1\hookrightarrow \D(\G_m)_{R/p^{2m}})\otimes_R R'
\ar[r]^-{f}_-{\sim} & (V(\G'_{2m})\hookrightarrow
\t{Lie}(E(\G'_{2m}))) }$$ We use Messing's Theorem (Theorem
\ref{messing}) to prove the above claim concerning the equivalence
of our deformation problem and the construction of of
\eqref{diag4.1}. Suppose we can construct such an admissible
filtration as in \eqref{diag4.1}. By Theorem \ref{messing}, we
have a \bt group $\G_{2m}$ over $R/p^{2m}$ which is a deformation
of $\G_m$ such that
$$(\t{Lie}(\G_{2m})\inj \t{Lie}(E(\G_{2m})) \simeq (\t{Fil}^1 \inj \D(\G_m)_{R/p^{2m}} ).$$
Thus $f$ induces an $R'$-isomorphism of \bt groups $\tilde f:
\G_{2m} \otimes_R R' \simeq \G'_{2m}$ which lifts the
Galois-equivariant isomorphism $\tilde f_0: \G_{m} \otimes_R R'
\simeq \G'_{m}$. Now it suffices to check that $\tilde f$ is
Galois equivariant. By Lemma \ref{check}, we need to show that
$\tilde f ^{-1}\circ \tilde f_\sigma =\t{Id}$. Since $f$, $f^1$
are Galois equivariant,  $f^{-1}\circ f_\sigma $ and
${(f^1)}^{-1}\circ f^1_\sigma$ are identity maps. By Theorem
\ref{messing} and the fact that the diagram \eqref{diag4.1} lifts
the diagram \eqref{diag4.2}, it is clear that $\tilde f ^{-1}\circ
\tilde f_\sigma =\t{Id}$. This settles one direction of the claim.
The proof of the other direction is just reversing the above
argument.

Now let us construct an  admissible filtration as in the diagram
\eqref{diag4.1}.
  Recall that for each $n$,  we have an exact sequence of $\Gamma$-descent
  modules
  \begin{equation}\tag{$S'_n$}
0 \rightarrow N'_n \rightarrow M'_n \rightarrow L'_n\rightarrow 0
\end{equation}
where $ ({N'_{n}\hookrightarrow M'_{n}})$ is just
 $({V(\G'_{n})\hookrightarrow \t{Lie}(E(\G'_{n}))})$.
 Over  $R/p^m$, by Remark \ref{remark},  we do have an isomorphism
of  effective descents of $(N'_m\hookrightarrow M'_m)$,
$$ {(N_m\hookrightarrow M_m) \simeq (V(\G_m)\hookrightarrow \t{Lie} (E(\G_m))}$$
Now choose a deformation $\tilde \G $ over $R/p^{2m} $ for $\G_m$
(such deformation always exists,  due to Theorem
\ref{Grothendieck}). Let $\D(\tilde \G )_{R/p^{2m}}=
\t{Lie}(E(\tilde \G ))$. From the functorial construction of the
crystal $\D(\tilde \G)_{R/p^{2m}}$, we have that $\D(\tilde \G
)_{R/p^{2m}} $ is a finite free $R/p^{2m }$-module and $\D(\tilde
\G)_{R/p^{2m}}\otimes_{R/p^{2m}}R/p^m\simeq \t{Lie}(E(\G_m))$. We
claim that $\D(\tilde \G)_{R/p^{2m}}$ is an effective descent for
$M'_{2m}=\t{Lie}(E(\G'_{2m}))$.
 Using  Lemma \ref{4.2} in our case, let
 $$(S_0\inj S) = (\t{Spec} {(R' /p^m)}\inj \t{Spec}(R'/ p^{2m}
 )),$$
 $$\G_0= \G_m\otimes_{R/p^m} {R'/p^m}, \ \H_0 = \G'_m ,\ \G = \tilde \G
 \otimes_{R/p^{2m}} R'/p^{2m} , \ \H = \G'_{2m}.$$
  We will prove
 the claim by showing that $E(\tilde \G)$ is an effective descent of $E(\G'_{2m})$.
 Since $\G_m $ is an effective descent of $\G'_m$, there exists an
 $\Gamma$-equivariant morphism $$f_0: \G_m\otimes_{R/p^m} {R'/p^m} \ito
 \G'_m.$$ Corollary  \ref{4.2.5} (2) shows that there exists a unique
 isomorphism $$E(f_0): E(\tilde \G)\otimes_{R/p^{2m}} {R'/p^{2m}} \ito E(\G'_{2m})$$ which lifts $f_0$.
Now it suffices to check $E(f_0)$ is $\Gamma$-equivariant. By
Lemma \ref{check}, we need to check $E(f_0)^{-1}\circ
E(f_0)_\sigma= \t{Id}$. But  this is clear by Corollary
\ref{4.2.5} (1) and the fact that ${f_0}^{-1}\circ
{(f_0)}_\sigma=\t{Id} $.

As  finite free $R/p^{2m}$-modules, $\D(\tilde \G )_{R/p^{2m}}$
has the same rank as $M_{2m}$. Now select an $R/p^{2m}$-module
isomorphism $f$ between
 $\D(\tilde \G )_{R/p^{2m}}$ and $M_{2m}$  such that  $f$ is the lift of isomorphism
$\D(\tilde \G )_{R/p^{2m}}\otimes_{R/p^{2m}}R/p^m\simeq M_m$.
There exists a finite free submodule $\t{Fil}^1$ of $\D(\tilde \G
)_{R/p^{2m}}$ such that the following diagram commutes
 \begin{equation}\label{diag4.3}
\begin{split}
 \xymatrix {\D(\tilde \G )_{R/p^{2m}}\ar[r]^-f_-{\sim}& M_{2m}\\
 \t{Fil}^1\ar@{^{(}->}[u] \ar[r]^{f^1}_{\sim}& N_{2m}\ar@{^{(}->}[u]
 }
\end{split}
 \end{equation}
 where $f^1$ is the restriction of $f$ on $\t{Fil}^1$.

 The map ${\t{Fil}^1\hookrightarrow \D(\tilde \G)_{R/p^{2m}}}$ is an admissible filtration because $L_{2m}\simeq M_{2m}/N_{2m}$
is a finite free $R/p^{2m}$-module.  By Lemma \ref{4.3} and the
diagram \eqref{diag4.3}, this admissible filtration is just what
we need in the diagram \eqref{diag4.1}. Therefore, we are done in
the case $s=2$. Note that above construction  gives an isomorphism
 $${(V(\G_{2m})\hookrightarrow \t{Lie}(E(\G_{2m})))\simeq (N_{2m}\hookrightarrow M_{2m})}$$
 as effective descents of  $ {N'_{2m}\hookrightarrow M'_{2m}}$,
 so we can continue our steps for  $s=3$ and then for any positive integers.
\end{proof}
Finally, by Lemma \ref{4.4}, defining  $\G =  \varinjlim_{s} \G_s
$, we see that  $\G$ is an effective descent for $\G'$. Thus,  we
have proved Theorem \ref{th1}.

\section{Extensions of Generic Fibres for Finite Flat Group Schemes}
\subsection{Proof of the Main Theorem}
In this Chapter, we assume that $R$ is a  complete discrete
valuation ring of mixed characteristic $(0, p)$ with perfect
residue field. We will generalize Tate's isogeny theorem \cite{ta}
and Raynaud' s results (Proposition 2.3.1 in \cite{ray}) to finite
level. First let us show how Theorem \ref{th2} and Proposition
\ref{th3} imply the Main Theorem. Then we will prepare to prove
Theorem \ref{th2} and Proposition \ref{th3}.
\begin{proof}[Proof of the Main Theorem] Applying Theorem \ref{th2} and
Proposition \ref{th3} to $R'$ and $G'=G \otimes_K K'$, we get
constants $c_1$ and $c_2$ which depend on $e(K'/\Q_p)$ and the
height of $G$ (note that heights of $G$ and $G'$ are the same).
Set $c=c_1+ c_2 +c_0\geq 2$. Suppose that $G[p^c]$ can be extended
to a finite flat group scheme $\G_{c}$. By Proposition \ref{th3},
$G[p^{c_1+ c_0 }]$ can be extended to a truncated \bt group $
\G_{c_1+c_0} $. Note that we have  $f: \G_{c_1+c_0} \otimes_R K'
\simeq \G'[p^{c_1 +c_0}]\otimes_{R'} K'$ Galois-equivariantly. By
Theorem \ref{th2}, $p^{c_1}f $ can be extended an $R'$-group
scheme morphism $F'$. Note that $F'$ kills the generic fiber of
$\G_{c_1+c_0} [p^{c_1}]\otimes_R R'$. Since $\G_{c_1+c_0}
[p^{c_1}] $ is flat,  $F'$ kills $\G_{c_1+c_0} [p^{c_1}]\otimes_R
R' $. Also,  it is easy to see that $F' $ factors though
$\G'[p^{c_0} ]$. Thus, $F'$ induces a morphism $F: \G_{c_1+c_0}
[p^{c_0}]\otimes_R R' \to \G' [p^{c_0}] $, and
$F\otimes_{R'}K'=f|_{\G_{c_1+c_0} [p^{c_0}]}$ is an isomorphism.

Since $f$ is Galois-equivariant and all group schemes here are
$R'$-flat,  any morphism extending $f$ has to be unique.
Therefore, $F$ is Galois-equivariant. Using the same argument for
$f^{-1}$, we see that there exists a morphism $$F^{-1}:
\G'[p^{c_0}]\to \G_{c_1+c_0} [p^{c_0}]\otimes_R R'.$$  Thus,
$\G_{c_1+c_0} [p^{c_0}] $ is an effective
 descent of $\G'[p^{c_0}]$. The Main Theorem then follows by the Descent Theorem (Theorem \ref{th1}).
 \end{proof}

\subsection{Dimension of truncated \bt groups} We will basically follow the similar ideas of Tate and Raynaud to prove our
theorem. So we first give an `` ad hoc" definition of the dimension for truncated \bt groups and then show
  that the dimension of truncated \bt groups can read off from the generic fiber if the level is big enough. Just like the case of \bt groups, this fact is
  crucial to extend the morphisms of generic fibers. Since
we will discuss the extension of generic fibers of finite flat
group schemes (truncated \bt group) defined over the same base
ring,  from now on, we fix our base ring $R$ which is a complete
discrete valuation ring of mixed characteristic (0, $p$) with
 perfect residue field and we assume that the absolute ramification index $e=e(K/\Q_p)$ is at least $p-1$.
We made such assumption because by Raynaud's Theorem in \cite{ray}, if $e< p-1$,
Theorem \ref{th2} and Theorem \ref{th3} automatically works for
$c_1=c_2=0$.   If $G$ is a scheme over $R$,
let $G_K$ denote the generic fiber of $G$ and same notations applies for the morphism of $R$-schemes.  For a finite flat group scheme $G$ over $R$, let
$\t{disc}(G)$ denote the discriminant ideal of $G$ and $|G|$ denote the order of $G$. Finally,
let $f: G \to H$ be a morphism of $R$-group schemes,  we call $f$ is a
\emph{generic isomorphism} if $f_K = f \otimes_R K$ is an $K$-isomorphism.

Now let us  introduce an invariant of truncated \bt groups, the dimension.
Given a truncated \bt groups $G$ over $R$, let $G_k= G \otimes_R k $ be its
closed fiber, $G_k^{(p)}=G_k\otimes_{k, \phi} k $,  where $\phi$ is
Frobenius. There exists a canonical $k$-group scheme morphism (relative Frobenious) $F: G_k\to G_k^{(p)}$ which kernel is
finite flat $k$-group scheme and $\t{Ker}(F) \inj  G[p]$ is a closed immersion.

\begin{definition}\label{def5.1}
 The dimension of $G$ is $\log_p |\t{Ker}(F)|$. We write  $\t{dim}(G)$ or $\t{d}(G)$ to denote it.
\end {definition}
Let us ``justify" Definition \ref{def5.1} by checking
compatibility with the definition of dimension for \bt groups. By
Grothendieck's Theorem in \cite{luc}, for any truncated \bt group
$G$ over $R$ of level $n$, there exists a \bt group $\G$ such that
$G= \G[p^n]$. We also have $\t{dim}(\G)= \log_p |\t{Ker}(F)|$
where $F: \G_k \to \G_k^{(p)}$ is the relative Frobenius on
$\G_k$,  so $\t{dim}(G) =\t{dim}(\G)$.

The following proposition show that if level $n$ is big enough,
the generic fiber of truncated \bt group will decide the dimension. For any $x$ is a real number, recall that
$[x]=\t{Max}\{m| m \t{ is an interger such that } m\leq x \}$ and $\lceil x\rceil= \t{Min}\{m| m \t{ is an
 integer such that } m\geq x \} $.
 \begin{prop}\label{dim}
Let $G$ and $H$ be  two truncated \bt group of level $n$  over $R$
which have the same generic fiber.  There exists $r_2\geq 2$
depending only on $e$ and height $h$ such that if $ n\geq r_2$  then
$\t{d}(G)=\t{d}(H)$.
 \end{prop}
 \begin{proof} It suffices to prove the proposition on the  strictly henselization of $R$. So we can assume that $k$ is
algebraically closed. From Grothendieck's Theorem in \cite{luc},
there exist \bt groups $\G$ and $\H$ such that $G= \G[p^n]$ and
$H=\H[ p^n ]$. In \cite{ray}, Raynaud proved that as one dimensional
Galois module $\t{det}(T_p({\G}))\simeq \epsilon ^d  $, where
$\epsilon $ is the cyclotomic character, $T_p (\G) $ is the Tate
module of $\G$ and $d=\t{dim}( \G)$. Since generic fibers of $G$ and
$H$ are the same, we have
$$\epsilon^{\t{d}(G)}=\t{det}((T_p({\G}))\equiv \t{det}(T_p(\H))=\epsilon^{\t{d}(H)}\mod(p^n)$$

Let $\lambda $ be the maximal integer such that $\mathbb Q_p(\zeta_{p^\lambda} )\subseteq K$, where $\zeta_{p^\lambda}$ is
$p^\lambda$-th root of unity.
Then  there exists a $\sigma \in \t{Gal}(\overline K/K)$ such that
$\epsilon(\sigma )= 1+p^\lambda u$, where $u$ is a unit in $\mathbb Z_p$.
Consider  the action of  $\sigma $ on the congruence relation above, we get
$$
(1+p^\lambda u )^{\t{d}(G)} \equiv (1+p^\lambda u)^{\t{d}(H)} \mod(p^n).
$$
  That is,   $\t{d}(G)\equiv \t{d}(H)\mod (p^{(n-\lambda)})$.
Note that $e(\Q_p(\zeta_{p^\lambda}):\Q_p)\leq
e$, we have $\lambda \leq [\log_p(\frac{e}{p-1})]+1. $ Set $r_2 =\lceil \log_ph\rceil + [\log_p(\frac{e}{p-1})]+1$.  we have
$$\t{d}(G)\equiv \t{d}(H)\mod (p^{(n-\lambda)})\equiv (p^{\lceil\log_p h\rceil } )$$
By proposition 6.2.8 in  \cite{bd3}, we know that $\t{d}(G),\  \t{d}(H) \leq h $, so  $\t{d}(G)=\t{d}(H)$.
\end{proof}

 From above Proposition and using the same argument of Tate (Proposition 6.2.12, \cite{bd3} or  \cite{ta}),
for truncated \bt  groups $G$ of level $n$,
 $\t{disc}(G)=p^{dnp^{nh}}$ for $n\geq r_2$, where  $d$ is the dimension
of $G$. Thus  we have following corollary:
 \begin{co}\label{genericiso}
Let  $f: G\rightarrow H $ be a morphism of truncated \bt groups of
level with $f_K$ an  isomorphism.
 If $n\geq r_2 $ then $f$ is an isomorphism.
\end{co}

The following lemma is a useful fact
on the scheme-theoretic closure of group schemes.

\begin{lemma}
Let $f: G\to H$ be a morphism of finite flat group schemes over
$R$,  $ G^{(1)}_K  $, $H_K^{(1)}$ be $K$-subgroup schemes of
$G_K$, $H_K$ respectively. Let $G^{(1)}$ and
$H^{(1)}$ be the scheme-theoretic closures of $G_K^{(1)}$,
$H_K^{(1)}$ in $G $ and $H$ respectively.
Suppose that $f\otimes_R K|_{G^{(1)}_K}$ factor through $H_K^{(1)}$.
Then $f|_{G^{(1)}}$ factors through $H^{(1)}$.
\end{lemma}
\begin{proof} Let $G'=G\times_{f, H, i }H^{(1)} $,
where $i : H^{(1)}\hookrightarrow H$ be the closed immersion. We see that $G'$ is a
closed subgroup scheme of $G$. It suffices
to show that the generic fiber of $G'$ contains $G_K^{(1)}$. Since $f_K|_{G^{(1)}_K}$ factor through $H_K^{(1)}$,
 this is clear.
\end{proof}

\begin{co}\label{5.6}
Notations as above, let  $G_i$, $H_i$ be scheme-theoretic closures
of $G_K[p^i ]$, $H_K[p^i]$ in $G$ , $H$ respectively, then $f|_{
G_i} $ factors through $H_i$.
\end{co}

\subsection{Proof of Proposition \ref{th3} and Theorem \ref{th2}}
Now we are ready to prove Proposition \ref{th3} and Theorem \ref{th2}
First we need to bound discriminants of finite flat group schemes
of type $(p,\cdots, p)$ in terms of heights and the absolute
ramification index $e(K/\Q_p)$.  

Let $G$ be a finite flat group scheme over $R$ killed by $p$. If
$|G|=p^h$, we define $h$ be the \emph{height} of $G$. Note that if
$\G$ is a \bt group (or a truncated \bt group) with height $h$.
Then $\G[p]$ will has height $h$ as a finite flat group scheme
killed by $p$. Recall that $\t{disc}(G)$ denote the discriminant
ideal of $G$.
\begin{lemma} With notations as above,   $\t{disc}(G)=(a)^{p^h}$, where $a\in R$ and
$a| p^h $.
\end{lemma}
\begin {proof}
Let  $0\to G_1 \to G_2 \to G_3 \to 0 $ be an exact sequence of
finite flat group schemes. By Lemma 6.2.9 \cite{bd3}, we have
\begin{equation}\label{disc}
 \t{disc} (G_2 ) =\t{disc}(G_1)^{|G_3|}\t{disc}
(G_3)^{|G_1|}
\end{equation}
Now passing to the connected-\'etale sequence $0 \to G^0 \to G \to
G^{\t{\'et}}\to 0$, since $\t{disc}(G^{\t{\'et}})=1$, we see that
$\t{disc}(G)= \t{disc}(G^0 )^{|G^{\t{\' et}}|} $. Thus, we reduce
the problem to the case that $G$ is connected.

  Now since $G$ is connected, by Theorem 3.1.1 in \cite{bbm}
 there exists an isogney $f$ between
 formal groups $\G$, $\G'$ such that
$$\xymatrix{0\ar[r]&G \ar[r] &\G\ar[r]^f&\G'\ar[r]& 0 }$$
is exact. By Corollary 4.3.10 in  \cite{bd3}, $\text{disc}
(G)=\text{N}_{\mathcal O_{\G}/f^*\mathcal O_{\G'}}(a)$, where
$\mathcal O_\G$ and $\mathcal O_{\G'}$ are formal Hopf algebras of
$\G$ and $\G'$ respectively;
 $(a)$ is the principal ideal in $\mathcal O_{\G'}$ for the annihilator of cokernel
$f^*(\Omega^n_{\G'})\rightarrow \Omega ^n_{\G}$ and $\text{N}_{\mathcal O_{\G}/f^*\mathcal O_{\G'}}(\cdot)$ is the
norm.

Let $d$ be the common dimension of $\G$ and $\G'$. By Theorem
3.2.1 in \cite{bd3}, we can choose basis of invariant
differentials $w'_1,\cdots,  w'_d$ and $w_1, \cdots,  w_d$ on
$\G'$ and $\G$ respectively. Since $f$ is an isogeny, for each
$i=1, \dots , d$, $f^*(w'_i)$ is an  invariant differential on
$\G$,  thus $f^*(w'_i)=\sum_{j=1}^d a_{ij}w_j$ for $a_{ij}\in R$.
 Define $a=\t{det}(a_{ij}) \in R $, so  $\text{disc} (G)=
\text{N}_{\mathcal O_{\G}/f^*\mathcal O_{\G'}}(a)= (a)^{p^h}$.

Now let us show that $a|p^h$. First let us reduce the problem to
the case that the generic fiber of $G$ is simple. In fact, suppose
$G'_{K}$ is a nontrivial subgroup of $G_K$, let $G'$ be the
scheme-theoretic closure of $G'_{K}$ in $G$. We get an exact
sequence of finite flat group schemes
$$ 0 \to G' \to G \to G'' \to 0. $$
Then by formula \eqref{disc}, we reduce the problem to the case
which the generic fiber $G$ is simple.

In general, if $G'\to G $ is a morphism of $R$-finite flat group
schemes that is an isomorphism on $K$-fibers,  we have
$\t{disc}(G)| \t{disc} (G')$. Thus,  if we fix the generic fiber
$G_K $ of $G$ it suffices to prove $a | p^h$ for the \emph{minimal
model} $G_{\t{min}}$ of $G_K$ in the sense of Raynaud. If $G_K$ is
simple, then by \cite{ray} such a minimal model is an $\bf
F$-group scheme of rank 1 for some finite field $\bf F$. Since
there is a connected model of $G$, clearly $G_{\t{min}}$ must be
connected.  Thus it suffices to prove that $a |p^h $ where $G$ is
a connected $\bf F$-group scheme of rank 1.

If we normalize the valuation by putting  $v(p)=e$, according to
Theorem 1.4.1 in  \cite{ray}
 the Hopf algebra $\mathcal O_G$ of a connected $\bf F$-group scheme has the following shape:
$$R\llbracket X_1,\cdots, X_h \rrbracket/(X_i^p-\delta_i X_{i+1}, i=1, \dots , h)$$
where $\delta_i\in R$ satisfies $1\leq v(\delta_i)\leq e$ for each
$i=1, \dots, h $.

Using Theorem 4.3.9 in \cite{bd3},  the \emph{different} of
$\mathcal O_G$ is the ideal of $R\llbracket X_1,\cdots,
X_h\rrbracket$ generated by
det$(\frac{\partial(X_i^p-\delta_iX_{i+1}) }{\partial X_i})$.
However,
$$
\t{det}(\frac{\partial(X_i^p-\delta_iX_{i+1}) }{\partial X_i})
=\begin{vmatrix}
pX_1^{p-1}& -\delta_1& \cdots & 0\\
0 & pX_2^{p-1}& \dots & 0\\
\cdots \\
-\delta_h& \cdots & \ &pX_h^{p-1}
\end{vmatrix}
=p^h\prod_{i=1}^{h}X_i^{p-1}-\prod_{i=1}^h\delta_i,
$$
so
 $$\t{disc}(G)=\t{N}_{\mathcal O_G/R}(p^h\prod_{i=1}^{h}X_i^{p-1}-\prod_{i=1}^h\delta_i).$$

Note that for any $X_i$, we have $X_i\cdot (X_1^{p-1}\cdots
X_h^{p-1})=\delta_1\cdots \delta_hX_i$ in  $\O_G $. Thus if we
select the $R$-basis of $\mathcal O_G$ given by  $ \{\prod_{i=1}^h
X_i^{n_i} , 0\leq n_i\leq p-1\}$ then for any $\prod_{i=1}^h
X_i^{n_i} $ with some $n_i>0$ we have the following equality in
$\O_G$:
$$(\prod_{i=1}^h X_i^{n_i})\cdot(p^h\prod_{i=1}^{h}X_i^{p-1}-\prod_{i=1}^h\delta_i) =(p^h-1)\prod_{i=1}^h\delta_i
 (\prod_{i=1}^h X_i^{n_i}) .$$
Thus,  $\t{disc}(G)=(-(p^h-1)^{p^h-1} (\delta_1 \cdots \delta_h
)^{p^h})$. Finally, since  $v(\delta_i )\leq e$ for each $i$
 we have  $\t{disc}(G)=(\delta_1 \cdots \delta_h )^{p^h }$ and $\delta_1 \cdots \delta_h | p^h $.
\end{proof}

\begin{proof}[Proof of Proposition \ref{th3}]
For each $1\leq i \leq n$, let $G_i $ be the scheme theoretic
closure of $G_K[p^i]$ in $G$, so  $G_i$ is a closed subgroup
scheme of $G_{i+1}$.

By Corollary \ref{5.6},  we see that $p^m$ induces a morphism
$G_{i+m}\to G_{i}$.   Furthermore,  $p^m $ induces morphisms
$\psi_{m,i, l}:G_{i+m+l}/G_{i+m}\to G_{i+l}/G_{i}$ which are
generic isomorphisms. Therefore, if $D_{i}$ is the affine algebra
of $G_{i+1}/G_{i}$, then for each $i$,
 $\psi_{1,1,i}:D_{i}\to D_{i+1}$ induces an
isomorphism $D_{i}\tensor K\simto D_{i+1}\tensor K$. If $A$
denotes the common affine algebra of the $(D_{i}\tensor K)$'s,
then we see by flatness that $D_{i}\to A$ is injective, so we can
identify $D_{i}$ with its image in $A$.  This  gives  a compatible
system of injections $D_{i}\inj D_{i+1}$. Therefore, we see that
$D_{i}$ is an increasing sequence of finitely generated
$R$-submodules of the finite \'etale $K$-algebra $A$. Since
$D_i\subset D_{i+1} $ inside $A$, we have
$\text{disc}(D_{i+1})|\t{disc}(D_i)\dots| \text {disc} (D_1)$. By
Lemma \ref{dim}, we see that $\t{disc}(D_i)=(a_i)^{p^h}$ and
$a_{i+1}| a_i| p^h $ for each $1\leq i \leq n$.
 Thus,  the number of total possible distinct
$D_i$ is at most $eh$.

Set $r_1= r_2(eh)+1$. If $n\geq r_1$, then there exists $i_0$ such
that
\begin{equation}\label{chain2}
D_{i_0}=D_{i_0+1}=\cdots =D_{i_0+s}\  \t{and}\   s \geq r_2.
\end{equation}
Let $\Gamma_{m}=G_{i_{0}+m}/G_{i_{0}}, m=1, \dots, s$. We claim
that $(\Gamma_{m})$ is a truncated {\bt} group of level $s$.

In fact, the closed immersions $G_i \hookrightarrow G_{i+1}$
induce closed immersions $i_m:\Gamma_m \hookrightarrow
\Gamma_{m+1}$. Since $(\Gamma_m) _K \simeq (G_K)[p^m]$,
$(\Gamma_m)$ is a generic truncated \bt group of level $s$. Thus,
all we have to check is that $i_{m}$ identifies $\Gamma_{m}$ with
$\Gamma_{m+1}[p^{m}]$ for all $m=0, \dots, s-1$. Consider the
diagram
\begin{equation}\label{E:final diag}
\begin{split}
\xymatrix{ \Gamma_{m+1}\ar@{=}[r] &
G_{i_{0}+m+1}/G_{i_{0}}\ar[r]^{p^{m}}\ar[d]^{\alpha} &
G_{i_{0}+m+1}/G_{i_{0}}\ar@{=}[r] & \Gamma_{m+1}\\
& G_{i_{0}+m+1}/G_{i_{0}+m}\ar[r]^{\beta} &
G_{i_{0}+1}/G_{i_{0}}\ar@{^{(}->}[u]^{\gamma}, & }
\end{split}
\end{equation}
where $\alpha$ is the canonical projection, $\beta$ is the map
induced by $p^{m}$, and $\gamma$ is the canonical closed
immersion.  Checking on the generic fiber, we see that
\eqref{E:final diag} commutes over $R$.  On the other hand, by our
choice of $i_{0}$ we see that $\beta$ is an isomorphism, and
therefore $\ker[p^{m}]=\ker\alpha$, which is nothing more than
$\Gamma_{m}$. This proves the claim that $(\Gamma_m)$ is a
truncated \bt group of level $s$.

Let $I$ be the set of indices $\{i| D_i=D_{i+1}=\cdots=D_{i+s_i},\
s_i\geq r_2  \}$,  $i''$ and $i'$  the maximal and minimal
elements of $I$ respectively. Let $\Gamma_m''= G_{i''+m }/G_{i''
},$ for $m=0,\dots, r_2$ and $\Gamma_m'= G_{i'+m }/G_{i' }$, for $
m=0,\dots, r_2$. We see that $(\Gamma_m'' )$ and $(\Gamma'_m)$ are
truncated \bt groups over $R$ of level $r_2$. Since $p^{i''-i'}$
induces a generic isomorphism $\Gamma_m''\to \Gamma_m'$,  from
Proposition  \ref{dim} we have $\t{ d} (\Gamma'')=\t{d}(\Gamma ' )
=d$, so  $D_{i''} = D_{i'}= p^{dp^h}$. Thus,  we have
\begin{equation}\label{chain}
D_{i'}=D_{i'+1}=\cdots = D_{i''}=\cdots =D_{i''+s_{i''}}.
\end{equation}

We claim that there are at least $n-r_1+r_2+1$ terms $D_i$ in
\eqref{chain}. In fact,
 the size of set of indices  $I '=\{i| i < i' \ \t{or } i > i'' +r_2\}  $
cannot be larger than $(eh-1)r_2 $ because otherwise,  since  we
have at most $eh$ many possible distinct $D_i$,  there would exist
a chain $D_i= D_{i+1}=\cdots = D_{i+ r_2}$ with $\{i,\dots
,i+r_2\}\subset I'$,  contradicting  the definition of $i '$ and
$i''$. Thus,
 there are at least $n-(eh-1)r_2=n-r_1 +r_2 +1 $ terms $D_i$ in \eqref{chain}.   Now  set $n_1= n-r_1+r_2$.
Let $G'_m= G_{i'+m}/G_{i'}$, for $  m=0, \dots , n_1$. Repeating
the proof above, we see that $(G'_m) $  is  a truncated \bt group
of level $n_1$ whose generic fiber is  $G_K[p^{n_1}]$.

Set $G'=G_{i'+ n_1}/G_{i'}$. We see that $p^{i'}$ induces a
generic isomorphism $\tilde g:G'\to G_{n_1} $. Let $g'$ be the
composite  of $\tilde g$ and the closed immersion
$G_{n_1}\hookrightarrow G$. We see that $g'_K $ factors though
$G_K[p^{n_1}]$,   and $g'_K : G' \to G_K[p^{n_1}]$ is an
isomorphism.

By Corollary \ref{5.6}, $ p^{r_1-r_2-i'} : G\to G $ factors
through $G_{i'+n_1}$, so composing with the natural projection
 $G_{i'+n_1}\to G_{i'+n_1}/G_{i'}=G'$ gives a morphism $g : G\to G'$ such that $g\otimes_R K =p^{r_1-r_2}$.

Finally, set $c_2=r_1-r_2$. The  above proof settles  the case
$n\geq r_1$. For the case $r_1\geq n \geq c_2$, let us return to
\eqref{chain2}. Since $r_2 \geq n-c_2$, it is easy to see that
there exists $i_0$ such that
$$D_{i_0}=D_{i_0+1}=\cdots =D_{i_0+s}\  \t{and}\   s \geq n-c_2. $$
 Let $\Gamma_m = G_{i_0 +m }/G_{i_0} $ for $m= 1, \dots
n-c_2$. We see that $(\Gamma_m)$ is a truncated \bt group of level
$n-c_2$. The rest of proof is the same as of the case $n\geq r_1$.
\end{proof}
\begin{proof}[Proof of Theorem \ref{th3}]
Set $c_1=r_1$. Let $M$ be the graph of $f$ in $G_K \times H_K$,
i.e., the scheme-theoretic image of $1\times f : G_K\to G_K\times
H_K$. It is obvious that $M$ is a truncated \bt group over $K$ of
level $n$. Using the natural projections of $G_K \times H_K $  to
$G_K$ and $H_K$, we have natural $K$-group scheme morphisms  $p_1
: M \to G_K $ and $p_2: M \to H_K$. Obviously, $p_1$ is an
isomorphism and $p_2\circ {p_1}^{-1}= f $.

For each $1\leq i \leq n $, let $E_i $ be the scheme-theoretic
closure of $M[p^i]$  in $G\times_R H $. From \S2, \cite{ray}, we
see that $E_i \inj E_{i+1}$ is a closed immersion and $E_i$ is a
closed subgroup scheme of $G\times_R H $.  We therefore have
morphisms of $R$-group schemes  $ \hat p_1  : E_n \to G$ and $\hat
p_2 : E_n \to H  $ which extend  morphisms $p_1:{ M}\to G_K$ and
$p_2: M\to H_K$ of generic fibers.
 Because $n\geq r_1 $, from Proposition \ref{th3} we see
that there exists a truncated \bt group $E'$ of level
$n_1=n-r_1+r_2$ and a morphism  $g': E'\to E_n $ such that $g'_K:
E' \to (E_n)_K[p^{n_1}]  $ is an isomorphism, so
  $\hat p_1 \circ g': E' \to G[p^{n_1}]  $ is an  isomorphism on generic fibers.
By Corollary \ref{genericiso}, since $n_1=n-r_1+r_2 \geq r_2$ we
see that $\hat p_1\circ g' $ is an isomorphism.  Finally, set $F=
{\hat  p_2 }\circ g' \circ (\hat p_1\circ g')^{-1}\circ p^{c_2} $.
We are done.
\end {proof}

\section{Application to Abelian Varieties }
\subsection{Proof of Theorem \ref{th1.2}} In this  chapter, we
assume that $R$ is a complete discrete valuation ring of mixed
characteristic $(0, p)$ with perfect residue field, and $K$ is the
fraction field.
 Let $A$ be an abelian variety defined over $K $
with dimension $g$. Suppose $A$  has potentially good reduction.
For example,
  $A$ can be taken to be an abelian variety with
 complex multiplication (\S 3,  \cite{serre}). Theorem \ref{th1.2} follows from the Main Theorem and the following lemma.

\begin{lemma} Suppose $A $ is an abelian variety over $K$ of dimension $g$ with potentially good reduction. There exists a constant $c_3$
which depends only  on $g$ such that  there exists a finite Galois
extension $K'$ of $K$  with $[K': K ]\leq c_3$ and  $A$ acquires
good reduction over $K'$.
\end{lemma}
\begin{proof} Choose a prime  $l \in \{ 3, 5\}$ with $l\not = p$.  Let $\rho_l: \gal(K^{\t{s}}/K)\to \t{Aut}{(T_l(A))}$
be the Galois representation. From \S 1, \cite{serre}, we know
that $A $ has potentially good reduction over $K$ if and only if
$\t{Im}(\rho_l)$ in $\t{Aut}{(T_l(A))}$ is finite.
  Since
$\t{Id}+l\t{Mat}_{2g}(\Z_l)$ has no torsion as multiplicative
group,
 we see that $\t{Id}+l\t{Mat}_{2g}(\Z_l)$ meets the finite image of $\rho_l $ trivially. Thus
$\bar \rho_l: \gal(K^{\t{s}}/K)\to \t{Mat}_{2g}{(\Z/l\Z)}$ is
injective. Therefore, the size of $\t{Im}({\rho}_l)$ is bounded by
the size of $\t{Mat}_{2g}{(\Z/l\Z)}$ which is $l^{4g^2}$. Let $K'$
be the fixed field by Ker($\rho_l$), so $A$ acquires good
reduction over $K'$ with $[K':K]=l^{4g^2}$. Finally, we can put
$c_3=\t{Max}\{3^{4g^2}, 5^{4g^2} $\}.
\end{proof}
\bibliographystyle{plain}

\bibliography{biblio}

\end{document}